\documentclass[ejs]{imsart}

\RequirePackage[OT1]{fontenc}
\RequirePackage{amsthm,amsmath}
\RequirePackage{natbib}
\RequirePackage[colorlinks,citecolor=blue,urlcolor=blue]{hyperref}
\usepackage{gensymb}
\usepackage{array,epsfig,fancyheadings,rotating}
\usepackage{sectsty, secdot}
\usepackage{colortbl,booktabs}
\usepackage{geometry}
\usepackage{amssymb}
\usepackage{amsfonts}
\usepackage{multirow}
\usepackage{bm, enumerate,xcolor}
\usepackage{epstopdf}
\usepackage{float,subcaption}
\usepackage{algorithm}
\usepackage{algorithmic}
\usepackage{url}

\allowdisplaybreaks

\newtheorem{theorem}{Theorem}
\newtheorem{lemma}{Lemma}
\newtheorem{corollary}{Corollary}

\theoremstyle{definition}

\pubyear{2020}
\volume{0}
\issue{0}
\firstpage{1}
\lastpage{34}

\startlocaldefs
\numberwithin{equation}{section}
\theoremstyle{plain}

\endlocaldefs

\begin{document}

\begin{frontmatter}
\title{Sparse Density Estimation with Measurement Errors}
\runtitle{Sparse Density Estimation with Measurement Errors}
%\thankstext{T1}{Footnote to the title with the `thankstext' command.}

\begin{aug}
\author{ \snm{Xiaowei Yang$^a$,}\ead[label=e0]{yxw8290@163.com}}
\author{ \snm{Huiming Zhang$^b$\thanksref{t2},}\ead[label=e1]{cute@pku.edu.cn}}
\author{\snm{Haoyu Wei$^c$}\ead[label=e2]{cute@pku.edu.cn}}
\and\author{\snm{Shouzheng Zhang$^d$}\ead[label=e3]{carol\underline{~}z@berkeley.edu}}

\address{a. School of Mathematics and Statistics, Chaohu University, Hefei, Anhui, China.\\
b.School of Mathematical Sciences, c.Guanghua School of Management, Peking University, Beijing, China.\\
 d. Department of Statistics, University of California, Berkeley, CA, USA.
\printead{e0,e1,e2,e3}}

%\author{\fnms{Third} \snm{Author}
%\ead[label=e3]{third@somewhere.com}
%\ead[label=u1,url]{www.foo.com}}
%
%\address{Address of the Third author\\
%usually few lines long\\
%usually few lines long\\
%\printead{e3}\\
%\printead{u1}}

%\thankstext{t1}{Some comment}
\thankstext{t2}{Corresponding author. The first three authors Xiaowei Yang, Huiming Zhang, and Haoyu Wei are co-first authors who contribute equally to this work. }
%\thankstext{t3}{Second supporter of the project}
%\thankstext{t1}{Some comment}
%\thankstext{t2}{First supporter of the project}
%\thankstext{t3}{Second supporter of the project}
\runauthor{Yang et al.}

\affiliation{Chaohu University, Peking University and University of California, Berkeley}

\end{aug}

\begin{abstract}
This paper aims to build an estimate of an unknown density of the data with measurement error as a linear combination of functions from a dictionary.
Inspired by the penalization approach, we propose the weighted Elastic-net penalized minimal $\ell_2$-distance method for sparse coefficients estimation, where the adaptive weights come from sharp concentration inequalities. The optimal weighted tuning parameters are obtained by the first-order conditions holding with a high probability. Under local coherence or minimal eigenvalue assumptions, non-asymptotical oracle inequalities are derived. These theoretical results are transposed to obtain the support recovery with a high probability. Then, some numerical experiments for discrete and continuous distributions confirm the significant improvement obtained by our procedure when compared with other conventional approaches. Finally, the application is performed in a meteorology data set. It shows that our method has potency and superiority of detecting the shape of multi-mode density compared with other conventional approaches.

\end{abstract}


\tableofcontents
\end{frontmatter}
%\newpage
\section{Introduction}\label{section 1}

Over the years, the mixture models have been extensively applied to model unknown distributional shapes in astronomy, biology, economics, and genomics see \cite{McLachlan19} and references therein. The distributions of real data involving potential complex variables often show multi-mode and heterogeneity. Due to the flexibility, it also often appear in various distribution-based statistical techniques, such as cluster analysis, discriminant analysis, survival analysis, empirical Bayesian inference. Flexible mixture models can naturally represent how the data are generated as mathematical artifacts. There are theoretical results showing that the mixture can approximate any density in the Euclidean space well, and the number of the mixture can also be finite (for example, a mixture of several Gaussian distributions). Although the mixture model is inherently attractive to the statistical modeling, it is a well-known difficult to infer, see \cite{Balakrishnanl17}. From the computational aspect, the optimization problems of mixture models are non-convex. Although existing computational methods, such as EM and various MCMC algorithms, are capable of making the mixture model fit the data relatively easily. It should be emphasized that the mixture problems are essentially difficult to be unrecognizable, and the number of components (says, the order selection) is hard to determine, see \cite{Chen2008}. There is a large amount of literature on its approximation theory and various methods have been proposed to estimate the components, see \cite{Das08} and references therein.

Nonparametric and combinatorial density estimation method were well studied in [\cite{Lugosi2001}, \cite{Biau2005}, as well as \cite{Martin09}]. These can be used to consistently estimate the amount of the components of the mixture when the components have known functional forms. When the number of candidate components is large, however, the non-parametric method becomes computationally infeasible. Fortunately, the advance of high-dimensional inference would compensate for this gap and guarantee the correct identification of the mixture components with a probability attending to one. With the advancement of technology, high-dimensional problems have been being applied to the forefront of statistical researches, and high-dimensional inference method has been applied to the infinite mixture models with a sparse mixture of $p \to \infty$ components, which is an interesting and challenging problem, see \cite{Bunea10} and \cite{Bertin11}. We propose an improvement of the sparse estimation strategy proposed in \cite{Bunea10}, in which Bunea et al. propose a $\ell_1$-type penalty to obtain a sparse estimate (SPADES), while we add a $\ell_2$-type penalty and extend the oracle-inequality results to our new estimator.

%However, difficult computing challenges have emerged in high-dimensional data, since it needs fast and flexible inference procedures. Calculations are quite cumbersome in high-dimensional problems, and the inference of the mixing distribution is still a challenging problem, see \cite{Martin09}.

In the real data, we often encounter the situation that the i.i.d. samples ${X_i} = {Z_i} + {\varepsilon _i}$ are contained by some zero-mean measurement error $\{ {\varepsilon _i}\}_{i = 1}^n$, see \cite{Chen1998},
\cite{Hall08}, \cite{Meister2006}, \cite{Cheng1999}. For density estimation of $\{X_{i}\}_{i=1}^{n}$, if there exists an orthogonal basis of functions, the estimation method is quite easy. In the measurement-error setting, however, finding an orthogonal based density function is not easy, see \cite{Schennach2013}. \cite{Schennach2013} suggests the assumption that the conditional distribution function of ${X_i}$ given ${Z_i}$ is known. This condition is somewhat strong since the most conditional distribution is hard to get the explicit formula (except the Gaussian distribution). To address this predicament, particularly with nonorthogonal base functions, the SPADES model is attractive and makes the situation easier to deal with. Based on the SPADES method, our approach is an Elastic-net calibration approach which is {simpler} and {more} interpretable than the conditional inference procedure proposed by \cite{Schennach2013}. In this paper, we proposed the corrected loss function to debase the measurement error, and this is motivated by \cite{Nakamura90}. We derive the honest variable selection consistency based on weighted $\ell_1$+ $\ell_2$ penalty, while some theoretical results of SPADES only contain the situation of the equal weights setting which is not plausible in sense of adaptive (data-dependent) penalized estimation. Moreover, we perform the Poisson mixture model to approximate the complex discrete distribution in the simulation part, while existing papers only emphasize the performance of continuous distribution models. Note that the multivariate kernel density estimator can only deal with continuous distribution and it requires a multivariate bandwidths section, while our method is dimensional free (the number of the required tuning parameters are only two).

This paper is presented as follows. Section {\ref{section 2}} introduces the density estimator which can deal with measurement errors. In this section, we introduce data-dependent weights for Lasso penalty, and the weights are derived by the event of KKT conditions holding with a high probability. In Section {\ref{section 3}}, we give a condition that can accurately estimate the weights of the mixture, with a probability tending to 1. We show that, in an increasing dimensional mixture model under the local coherence assumption, if the tuning parameter is higher than the noise level, the recovery of the mixture component can hold with a high probability. In Section {\ref{section 4}}, we study the performance of our approach on artificial data generating from mixed Gaussian or Poisson distributions compared with other conventional methods, which indeed shows the improvement by employing our procedure. Besides, the simulation also demonstrates that our method is better than the traditional EM algorithm even under a low dimensional model. Considering the multi-modal density aspect of the meteorology dataset, our proposed estimator has a stronger ability in detecting multiple modes for the underlying distribution, comparing with other methods such as SPADES or un-weighted Elastic-net estimator. Section {\ref{section 5}} is the summary and the proof of theoretical results is delivered in the Appendix.
%illustrates our methodology through numerical simulations and real data sets.
%Real data applications are presented in section 6.
%Section \ref{sec:summary} concludes this paper with some discussions.
%All proofs are delivered to the Appendix.

\section{Density Estimation}\label{section 2}

%\textbf{Notations}:
% For the sake of brevity, we summarize some of the basic symbols that will be used in the rest of this paper.  We considered the case of nonrandom covariates. In what follows, for any vector $\bm x=(x_0,x_1,\ldots,x_p)^T \in \mathbb{R}^{p}$, let $\|\bm x\| :=  (\sum_{j = 0}^p {x_j^2})^{1/2}$ denote the Euclidean norm and let $\|\bm x\| := \sum_{j = 0}^p {|{x _j}|}$ denote the $l_1$ norm. Let $\|\bm x\|_\infty:=\mathop {\max }_{0 \le j \le p} |{x _j}|$ denote the $l_\infty$-norm. ``Constant" refer to finite, positive, and non-random numbers.
\subsection{Mixture models}\label{sec:glm}
Suppose that $\{Z_{i}\}_{i=1}^{n}$ are independent random variables with a common unknown density $h \in \mathbb{R}^{d}$. However, the observations are contaminated with measurement errors $\{ {\varepsilon _i}\}_{i = 1}^n$ as latent variables, the observed data are actually ${X_i} = {Z_i} + {\varepsilon _i}$. Let $\{h_{j}\}_{j=1}^{W}$ be a series of density functions (such as Gaussian or Poisson), and $\{h_{j}\}_{j=1}^{W}$ are also called basis functions. Assume that the estimator of $h$ belongs to the linear combination of $\{h_{j}\}_{j=1}^{W}$. The $W:=W_n$ is a function of $n$, which is of particular intrigueing for us, since there may be $W \gg n$ (the high-dimensional setting). Let $\beta^*:=(\beta_{1}^*,\cdots,\beta_{W}^*)\in \mathbb{R}^{W}$ be the unknown true parameter. Assume that
\begin{itemize}
\item [\textbullet] (H.1): the $h:=h_{\beta^*}$ is defined as
\begin{equation}\label{eq:1}
Z \sim h(z):=h_{\beta^*}(z)=\sum_{j=1}^{W}\beta_{j}^*h_{j}(z),~~\text{with}~\sum\limits_{j = 1}^W {\beta _j^*}= 1.
\end{equation}
\end{itemize}
If the base is orthogonal and there are no measurement errors, a perfectly natural method is to estimate $h$ by an orthogonal series of estimators in the form of $h_{\tilde{\beta}}$, where $\tilde{\beta}$ has the coordinates $\tilde{\beta}_j=\frac{1}{n}\sum_{i=1}^{n}h_{j}(X_{i})$. However, this estimator depends on the choice of $W$, and a data-driven selection of $W$ or the threshold needs to be adaptive. This method, as well as many other methods, can only be applied to $W \leq n$, nonetheless, we want to solve more general problems, the base functions $\{h_{j}\}_{j=1}^{W}$ are not necessarily orthogonal for instance. Here the $W$ is not necessarily less than $n$, and we want to achieve the best convergence.

Theorem 33.2 in \cite{Das08} states that any smooth density can be well approximated by a finite mixture of some continuous functions. However, Theorem 33.2 in \cite{Das08} does not make sure that how many components $W$ are required for the mixture. Thus the hypothesis of the increasing-dimensional $W$ is reasonable. For discrete distributions, there is also a similar mixture density approximation, see Remark of Theorem 33.2 in \cite{Das08}. Suppose that:

\begin{itemize}
\item [\textbullet] (H.2): the density of the observed data is the linear combination of a series of some new based density functions $\{{\tilde h}_{j}\}_{j=1}^{W}$  if we still use the original true parameter $\beta^*$:
\begin{equation}\label{eq:1g}
X \sim g(x):=g_{\beta^*}(x)=\sum_{j=1}^{W}\beta_{j}^*{\tilde h}_{j}(x).
\end{equation}
\end{itemize}
Note that the $\{\tilde h_{j}\}_{j=1}^{W}$ are not the true mixture density since we wrongly use the contaminated data to fit the unobserved data. The direct estimation of coefficients based on $\{\tilde{h}_{j}\}_{j=1}^{W}$ in equation \eqref{eq:1g} is nothing but imprecise.

\subsection{Estimation for sparse density with measurement errors}
This subsection aims to construct a sparse estimator for the density $h(z):=h_{\beta^*}(z)$ as a linear combination of known densities.

Recall the definition of the $L_{2}(\mathbb{R}^{d})$ norm $\|f\|=\left(\int_{\mathbb{R}^{d}}f^{2}(x)dx\right)^{\frac{1}{2}}$. For $f,g\in L_{2}(\mathbb{R}^{d})$, let the inner product be $<f,g>=\int_{\mathbb{R}^{d}}f(x)g(x)dx$. Note that if the density $h(z)$ belongs to $L_{2}(\mathbb{R}^{d})$ and assume that $\{X_{i}\}_{i=1}^{n}$ has the same distribution $X$, for any $f\in L_{2}$, we have $<f,h>=\int_{\mathbb{R}^{d}}f(x)h(x)dx=Ef(X)$. If $h(x)$ is the density function for a discrete distribution, the integral is replaced by summation, and we can define the inner product as $< f,h > : = \sum\nolimits_{k \in \mathbb{Z}^{d}} {f(k)h(k)}$.

Let us minimize the $\|h_{\beta}-h\|^{2}$ on $\beta \in \mathbb{R}^{W}$ to obtain the estimate of $h(z):=h_{\beta^*}(z)$, i.e. minimizing
\begin{align}\label{eq:2}
\|h_{\beta}-h\|^{2}&=\|h\|^{2}+\|h_{\beta}\|^{2}-2<h_{\beta},h>\nonumber\\
&=\|h\|^{2}+\|h_{\beta}\|^{2}-2Eh_{\beta}(Z) \propto -2Eh_{\beta}(Z)+\|h_{\beta}\|^{2}.
\end{align}
The \eqref{eq:2} implies that minimizing the $\|h_{\beta}-h\|^{2}$ for true observations $\{Z_{i}\}_{i=1}^{n}$ is equivalent to minimizing
\begin{align}\label{eq:appro}
-2Eh_{\beta}(Z)+\|h_{\beta}\|^{2}\approx  - \frac{2}{n}\sum\limits_{i = 1}^n {{h_\beta }({Z_i})}+\|h_{\beta}\|^{2}.
\end{align}
For the convenience analysis of the measurement error, suppose that
\begin{itemize}
\item [\textbullet] (H.3): the i.i.d. observations $\{X_{i}\}_{i=1}^{n}$ have following decomposition:
\[g_{\beta}(X_1)=\sum_{j=1}^{W}\beta_{j}{\tilde h}_{j}(X_1) = \sum\limits_{j = 1}^W {{\beta _j}} [{h_j} * {f_\varepsilon}](X_1)\approx \sum\limits_{j = 1}^W {{\beta _j}} [{h_j} + {e_j}](X_1)\]
where $\{{e_j}(x)\}_{j=1}^W$ are technically assumed to be some orthogonal error functions, and the given base functions $\{ h_{j}\}_{j=1}^{W}$ are orthogonal to the error functions $\{{e_j}(x)\}_{j=1}^W$. Moreover, we also assume that the error functions as the perturbation functions have the zero empirical average evaluated at the observed data $\frac{1}{n}\sum\limits_{i = 1}^n {{e_j}({X_i})}  \approx E{e_j}(X) = 0$ for any $j$.
\end{itemize}

The assumption (H.3) means that $\{{e_j}(x)\}_{j=1}^W$ is an instrumental function to deal with the misspecified base function $\{\tilde h_{j}\}_{j=1}^{W}$ in \eqref{eq:1g}. Here, we mimic the idea of instrumental variables in econometrics, so we suppose that there always exists such functions. For $g_{\beta}(X_1)$ we have the following approximation:
\[ - 2E{g_\beta }(X) \approx - 2E\sum\limits_{j = 1}^W {{\beta _j}} [{h_j} + {e_j}](X) =  - 2E{h_\beta }(X) \approx  - \frac{2}{n}\sum\limits_{i = 1}^n {{h_\beta }({X_i})} .\]

For the observations $\{X_{i}\}_{i=1}^{n}$ with the measurement errors $\{ {\varepsilon _i}\}_{i = 1}^n$, minimizing the $\|g_{\beta}-g\|^{2}$  is equivalent to minimizing $-2Eg_{\beta}(X)+\|g_{\beta}\|^{2}$. More specifically, we approximate $-2Eg_{\beta}(X)+\|g_{\beta}\|^{2}$ by the argument below:
\begin{align*}\nonumber
-2Eg_{\beta}(X)+\|g_{\beta}\|^{2}& \approx  - \frac{2}{n}\sum\limits_{i = 1}^n {{h_\beta }({X_i})}+\sum\limits_{1 \le i,j \le W} {{\beta _i}} {\beta _j}\int_{{\mathbb{R}^d}} {[{h_i}(z) + {e_i}(z)][{h_j}(z) + {e_j}(z)]dz} \\
&= -\frac{2}{n}\sum\limits_{i = 1}^n {{h_\beta }({X_i})}+\sum\limits_{1 \le i,j \le W} {{\beta _i}} {\beta _j}\int_{{\mathbb{R}^d}} {{h_i}(z){h_j}(z)dz}  + \sum\limits_{1 \le i \le W} {\beta _i^2} \int_{{\mathbb{R}^d}} {e_i^2(z)dz}\\
&\approx -\frac{2}{n}\sum\limits_{i = 1}^n {{h_\beta }({X_i})}+\|h_{\beta}\|^{2} + c \sum\limits_{1 \le i \le W} {\beta _i^2},
\end{align*}
where the equality stems from the orthogonality assumption for $\{{e_j}(x)\}_{j=1}^W$. It is easy to see that $c$ introduces some approximate information regarding the measurement errors, see \cite{Rosenbaum2010} for similar purpose. It is different from SPADES, since adjusting for the presence of measurement error is important for accurately describing the relationship between the true covariates and the outcome of interest.

It is plausible to assign more constrains for the candidate set of $\beta$ in the optimization, for example, the $\ell_1$ constrains $\|\beta\|_{1} \le a$ where $a$ is the tuning parameter. More adaptively, we prefer to use the weighted $\ell_1$ restriction $\sum_{j=1}^{W}\omega_{j}|\beta_{j}|\le \tilde{a}$ where the weights $\omega_{j}$'s are data-dependent will be specified later. From the discussion above, now we propose the following \emph{Corrected Sparse Density Estimator} (CSDE):
\begin{eqnarray}\label{eq:9}
\hat{\beta}:=\hat{\beta}(\omega_{1},\cdots,\omega_{W})=\arg\min_{\beta\in \mathbb{R}^{W}}\left\{-\frac{2}{n}\sum_{i=1}^{n}h_{\beta}(X_{i})+\|h_{\beta}\|^{2}+2\sum_{j=1}^{W}\omega_{j}|\beta_{j}|+c\sum_{j=1}^{W}\beta_{j}^{2}\right\}
\end{eqnarray}
where $c$ is the tuning parameter for $\ell_{2}$-penality. As we said before, $c$ also present the correction for adjusting the measurement errors in our observations.

For CSDE, if $\{h_{j}\}_{j=1}^{W}$ is an orthogonal system, it can be clearly seen that the CSDE estimator is consistent with the soft threshold estimator, and the explicit solution is ${{\hat \beta }_j} = \frac{{{{(1 - {\omega _j}/|{{\tilde \beta }_j}|)}_+ }{{\tilde \beta }_j}}}{{1 + c}}$, where $\tilde{\beta}_{j}=\frac{1}{n}\sum_{i=1}^{n}h_{j}(X_{i})$ and $x_{+}=\max(0,x)$. In this case, we can see that $\omega_{j}$ is the threshold of the $j$-th component of the simplest mean estimator $\tilde{\beta}=(\tilde{\beta}_{1},\cdots,\tilde{\beta}_{W})$.

From sub-differential of the convex optimization, the corresponding Karush-Kuhn-Tucker conditions (necessary and sufficient first-order condition) for minimizer \eqref{eq:9} is
\begin{lemma}[KKT conditions in short, Lemma 4.2 of \cite{Buhlmann2011}]\label{lemma 1}
 Let $k \in \{ 1,2, \cdots ,W\}$ and $c > 0$. Then, a necessary and sufficient condition for CSDE to be a solution of \eqref{eq:9} is
\begin{enumerate}
\item
${\hat \beta _k}: \ne 0$ if
$
\frac{1}{n}\sum_{i=1}^{n}h_{k}(X_{i})-\sum_{j=1}^{W}\hat{\beta}_{j}<h_{j},h_{k}>-c\hat{\beta}_{k}=w_{k}{\rm{sign}}(\hat{\beta}_{k}).$
\item
${\hat \beta _k}= 0$ if
$\left|\frac{1}{n}\sum_{i=1}^{n}h_{k}(X_{i})-\sum_{j=1}^{W}\hat{\beta}_{j}<h_{j},h_{k}>-c\hat{\beta}_{k}\right|\leq w_{k}.$
\end{enumerate}
\end{lemma}

Since all $\beta_{j}^{*}$ are non-negative, when doing minimization in equation (\ref{eq:9}), we have to put a non-negative restriction for optimizing \eqref{eq:9}.
%
%For example, if $\{h_{j}\}_{j=1}^{W}$ orthogonal system, it can be clearly seen that the CSDE estimator is consistent with the soft threshold estimator, and the explicit solution is ${{\hat \beta }_j} = \frac{{{{(1 - {\omega _j}/|{{\tilde \beta }_j}|)}_+ }{{\tilde \beta }_j}}}{{1 + c}}$, where $\tilde{\beta}_{j}=\frac{1}{n}\sum_{i=1}^{n}h_{j}(X_{i})$ and $x_{+}=\max(0,x)$. In this case, we can see that $\omega_{j}$ is the threshold of the $j$th component of the simple mean estimator $\tilde{\beta}=(\tilde{\beta}_{1},\cdots,\tilde{\beta}_{W})$.

In fact, we do prefer to adapt the weight lasso penalty as a convex adaptive $\ell_1$ penalization due to the computational feasibility and optimal first-order conditions. We require that the larger weights are assigned to the coefficients of unimportant covariates, while the smaller weights are accompanied by important covariates. So the weights represent the importance of the covariates. The larger (smaller) weights shrink to zero more easily (difficultly) than the un-weighted lasso, with appropriate or even optimal weights, which may lead to less bias and more efficient variable selection. The derivation of the weight will be given in Section {\ref{section 2.4}}.

\subsection{We cannot transform the mixture models to linear models!}

In this part, we will illustrate that in the mixture models, even without measurement error, \eqref{eq:1} can't be partially transformed into the linear model, namely
\begin{equation}\nonumber
Y={\rm{X}}^{T}\beta+\varepsilon,
\end{equation}
where $Y$ is the $n$-dimensional response variables, ${\rm{X}}$ is the $W\times n$-dimensional fixed design matrix, ${\beta}$ is a $W$-dimensional vector of model parameters,  the ${\varepsilon}$ is a $n\times1$-dimensional vector for random error terms with zero mean and finite variance. Consider the least square objective function $U({\beta})$ for estimating ${\beta}$,
\begin{align}\label{eq:Q}
U({\beta})&=(Y-{\rm{X}}^{T}\beta)^{T}(Y-{\rm{X}}^{T}\beta)=-2Y^{T}{\rm{X}}^{T}\beta+\beta^{T}{\rm{X}}{\rm{X}}^{T}\beta+Y^{T}Y.
\end{align}
Minimizing (\ref{eq:Q}) is equivalent to minimizing $U^{*}(\beta)$ in the following formula (\ref{eq:Qstare})
\begin{equation}\label{eq:Qstare}
U^{*}(\beta)=-2Y^{T}{\rm{X}}^{T}\beta+\beta^{T}{\rm{X}}{\rm{X}}^{T}\beta.
\end{equation}
Comparing the objective function \eqref{eq:Qstare} with \eqref{eq:appro}, it is easy to obtain
\begin{equation}\nonumber
Y=(\frac{1}{n},\frac{1}{n},\cdots,\frac{1}{n})^{T}, ~ \beta=(\beta_{1},\beta_{2},\cdots,\beta_{W})^{T}, ~
{\rm{X}}=\left(
\begin{array}{ccc}
h_{1}(X_{1}) & \cdots & h_{1}(X_{n}) \\
\vdots & \ddots & \vdots \\
h_{W}(X_{1}) & \cdots & h_{W}(X_{n})
\end{array}
\right).
\end{equation}
Substituting $Y$, ${\rm{X}}$ and $\beta$ into a linear regression model, we get
\begin{equation}\nonumber
\left(
\begin{array}{c}
\frac{1}{n}\\
\vdots \\
\frac{1}{n}
\end{array}
\right)_{n\times1}=\left(
\begin{array}{ccc}
h_{1}(X_{1}) & \cdots & h_{W}(X_{1}) \\
\vdots & \ddots & \vdots \\
h_{1}(X_{n}) & \cdots & h_{W}(X_{n})
\end{array}
\right)_{n\times W}\left(
\begin{array}{c}
\beta_{1} \\
\vdots \\
\beta_{W}
\end{array}
\right)_{W\times1}+\left(
\begin{array}{c}
\varepsilon_{1} \\
\vdots \\
\varepsilon_{n}
\end{array}
\right)_{n\times1}.
\end{equation}
Then,
\begin{equation}\label{eq:88}
\varepsilon_{i}=\frac{1}{n}-\sum_{j=1}^{W}\beta_{j}h_{j}(X_{i}), ~ i=1,2,\cdots,n.
\end{equation}
It can be seen from equation (\ref{eq:88}) that the value of $\varepsilon_{i}$ is no longer random if ${\rm{X}}$ was the fixed design matrix. Furthermore, even ${\rm{X}}$ is a random design, take the expectation on both sides of (\ref{eq:88}), and one can find that the left side is not equal to the right side, that is, $$E(\varepsilon_{i})= 0 =\frac{1}{n} - \sum\limits_{j = 1}^W {{\beta _j}} E{h_j}({X_i}).$$
It leads to additional requirement $\sum\limits_{j = 1}^W {{\beta _j}} E{h_j}({X_i}){\rm{ = }}\frac{1}{n} \to 0$ which is meaningless as $n\rightarrow\infty$, since all ${\beta _j}$ and ${h_j}$ are positive, this is a contradiction to $\sum\limits_{j = 1}^W {{\beta _j}} E{h_j}({X_i})>0$ for all $n$.

Both of the two situations above contradict the definition of the assumed linear regression model, hence we can't convert the estimation of \eqref{eq:1} into the estimation problem of linear models. Thus the existing oracle inequalities are not applicable anymore, we will propose the oracle inequalities later. However, we can transform the mixture models to corrected score Dantzig selector such as \cite{Belloni17}. Although \cite{Bertin11} studies the oracle inequalities for adaptive the Dantzig density estimation, their study does not contain the error-in-variables framework and the support recovery content.

\section{Sparse Mixture Density Estimation}\label{section 3}
In this section, we will present the oracle inequalities for estimators $\hat{\beta}$ and $h_{\hat{\beta}}$. The core of this section consists of 5 main results, corresponding to the oracle inequalities for estimated density (Theorems {\ref{theorem 1}} and {\ref{theorem 2}}, respectively), upper bounds on $\ell_1$-estimation error (Corollaries {\ref{corollary 1}} and {\ref{corollary 2}}, respectively), and support consistency (Theorem {\ref{theorem 3}}) as the byproduct of Corollary {\ref{corollary 2}}.

\subsection{Data-dependent weights}\label{section 2.4}
The weights $\omega_{j}$'s are chosen adequately such that the KKT condition for stochastic optimization problems has a high probability to be satisfied.

As mentioned before, the weights in (\ref{eq:9}) rely on the observed data since we calculate the weights which make sure the KKT conditions hold with a high probability. The weighted Lasso estimates could have less $\ell_1$ estimation error comparing with Lasso estimates, see also the simulation part. Next, the question we need to consider is what kind of data-dependent weights configuration can enable the KKT conditions to be satisfied with a high probability. The fundamental way to get data-dependent weights is to apply a concentration inequality for a weighted sum of independent random variables. Moreover, the weights should be a known function of data without any unknown parameter. There is a criterion that can help to obtain the weight grounded on Bernstein's concentration inequality in SPADES. Whereas, the convergence rate of the probability upper bounds of the summation of $n$ independent random variables deviated from its expected value for Bernstein's concentration inequality is $\exp \left( { - \frac{{{c_1}{t^2}}}{{{c_2 n} + {c_3}t}}} \right)$. Contrasting to the Bernstein's concentration inequality, the McDiarmid's inequality (also known as the bounded difference inequality which is used for obtaining the desired weights) has a faster convergence rate $\exp \left( { - \frac{{{c_1}{t^2}}}{n}} \right)$ in $t$.

\begin{lemma}\label{lemma 2}
Suppose $X_{1},\cdots,X_{n}$ are  independent random variables, all values belong to $A$. Let $f:A^{n}\rightarrow \mathbb{R}$ be a function and satisfy the bounded difference conditions
\begin{equation}\nonumber
\sup_{x_{1},\cdots,x_{n},x_{s}^{'} \in A }|f(x_{1},\cdots,x_{n}) - f(x_{1},\cdots,x_{s-1},x_{s}^{'},x_{s+1},\cdots,x_{n})|\leq C_{s},
\end{equation}
then for all $t>0$,
\begin{equation}\nonumber
P\left\{|f(X_{1},\cdots,X_{n})-Ef(X_{1},\cdots,X_{n})|\geq t\right\}\leq 2\exp\left\{-\frac{2t^{2}}{\sum_{s=1}^{n}C_{s}^{2}}\right\}.
\end{equation}
\end{lemma}
We define the KKT conditions of optimization evaluated at $\beta^*$ (it is from the sub-gradient of the optimization function evaluated at $\beta^*$) by the events below:
\begin{equation}\nonumber
{\cal{F}}_{k}(\omega_{k}):=\left\{\left|\frac{1}{n}\sum_{i=1}^{n}h_{k}(X_{i})-\sum_{j=1}^{W}\beta_{j}^*<h_{j},h_{k}>-c {\beta _k^*}\right|\leq\omega_{k}\right\},k=1,2,\cdots,W.
\end{equation}

Assume that
\begin{itemize}
\item [\textbullet] (H.4): $\exists L_k>0$ s.t. ${\left\| {{h_k}} \right\|_\infty } = \mathop {\max }\limits_{1 \le i \le n} |{h_k}({X_i})|\leq L_k$.
\end{itemize}
\begin{itemize}
\item [\textbullet] (H.5): $0<\mathop {\max }\limits_{1 \le j \le W} |\beta _j^*| \le B$.
\end{itemize}
The (H.4) and (H.5) is a common a assumption in sparse $\ell_1$ estimation, see \cite{Bunea10}, \cite{Zhang20}.

Hence we could check that the following event is verified by some difference conditions. Note that $Eh_{k}(X_{i})=\sum_{j=1}^{W}\beta_{j}<h_{j},h_{k}>$ (which is free of $X_{i}$), we have
\begin{align*}
\frac{1}{n}&\left|\sum_{i=1}^{n}h_{k}(X_{i})-\sum_{j=1}^{W}\beta_{j}^*<h_{j},h_{k}>-\left(\sum_{i\neq s}^{n}h_{k}(X_{i})-\sum_{j\neq s}^{W}\beta_{j}^*<h_{j},h_{k}>+h_{k}(X_{s}^{'})-Eh_{k}(X_{s}^{'})\right)\right|\\
&=\frac{1}{n}\left|h_{k}(X_{s})-h_{k}(X_{s}^{'})+Eh_{k}(X_{s}^{'})-Eh_{k}(X_{s})\right|\\
&\leq\frac{1}{n}(|h_{k}(X_{s})-h_{k}(X_{s}^{'})|+|Eh_{k}(X_{s}^{'})-Eh_{k}(X_{s})|)\leq\frac{4L_{k}}{n}.
\end{align*}
The last inequality above is due to $|h_{k}(X_i)-Eh_{k}(X_i)|\leq 2L_{k}$.

Next, we apply the McDiarmid's inequality on the event ${\cal{F}}_{k}^{c}(\omega_{k})$ by (H.5). Then
\begin{align*}
P({\cal{F}}_{k}^{c}(\omega_{k}))&=P\left\{\left|\frac{1}{n}\sum_{i=1}^{n}h_{k}(X_{i})-\sum_{j=1}^{W}\beta_{j}^*<h_{j},h_{k}>-c {\beta _k^*}\right|\geq \omega_{k}\right\}\\
&\leq P\left\{\left|\frac{1}{n}\sum_{i=1}^{n}h_{k}(X_{i})-Eh_{k}(X_{i})\right|+c {\beta _k^*}\geq \omega_{k}\right\}\\
(\text{by~(H.5)})~&\leq P\left\{\left|\frac{1}{n}\sum_{i=1}^{n}h_{k}(X_{i})-Eh_{k}(X_{i})\right|\geq \omega_{k}-cB\right\}\\
(\text{define}~\tilde\omega_{k}:=\omega_{k}-cB >0)~& \le 2\exp \left\{ { - \frac{{2\tilde \omega _k^2}}{{16L_k^2/n}}} \right\}=2\exp\left\{-\frac{n\tilde\omega_{k}^{2}}{8L_{k}^{2}}\right\}=:\frac{\delta}{W}, \quad 0<\delta<1.
\end{align*}
Considering the previous line,
\begin{equation}\label{eq:3.1}
\omega_{k}:=2\sqrt{2}L_{k}\sqrt{\frac{1}{n}\log \frac{2W}{\delta}}+cB=:2\sqrt{2}L_{k}v(\delta/2)+cB,
\end{equation}
where $v=v(\delta ): = \sqrt {\frac{1}{n}\log \frac{W}{\delta }}$.

The weight $\omega_{k}$ in our paper is different from \cite{Bunea10} that gives the un-shift version ($\check{\omega_{k}}=4L_{k}\sqrt{\frac{1}{n}\log \frac{W}{\delta/2}}$), due to the Elastic-net penalty. And we define a modified version of event of KKT condition
\begin{equation}\label{eq:kkt1}
{\cal{K}}_{k}(\omega_{k}):=\left\{\left|\frac{1}{n}\sum_{i=1}^{n}h_{k}(X_{i})-\sum_{j=1}^{W}\beta_{j}^*<h_{j},h_{k}>\right|\leq\tilde\omega_{k}\right\},k=1,2,\cdots,W
\end{equation}
holds with the probability $1-2\exp\left\{-\frac{n\tilde\omega_{k}^{2}}{8L_{k}^{2}}\right\}$ at least.

\subsection{Non-asymptotic oracle inequalities}
The oracle inequality connects the performance of an obtained estimator with the true parameter which is not available in practice, see \cite{Candes08} for more discussions. Introduced by \cite{Donoho1994}, oracle inequality is a powerful non-asymptotical and analytical tool that seeks to provide the distance from the obtained estimator and a true estimator.

For $\forall\beta\in \mathbb{R}^{W}$, let
\begin{equation}\nonumber
I(\beta)=\{j\in \{1,\cdots,W\}:\beta_{j}\neq0\}
\end{equation}
be the indices corresponding to the non-zero components of the vector $\beta$, i.e. the support in mathematical jargon. If there is no ambiguity, we would like to write $I(\beta^{*})$ as $I_{*}$ for simplicity. And
\begin{equation}\nonumber
W(\beta)=\sum_{j=1}^{W}I(\beta_{j}\neq0)
\end{equation}
is the number of its non-zero components, where $I(\cdot)$ represents the indicative function. Let
\begin{equation}\nonumber
\sigma_{j}^{2}=Var(h_{j}(X_{1})), \quad 1\leq j\leq W.
\end{equation}

Below, we will state the non-asymptotic oracle inequalities for $h_{\hat{\beta}}$ (with the high probability  $1-\delta(W,n)$ for any integer $W$ and $n$) which measures the $L_2$ distance between $h_{\hat{\beta}}$ and $h$.  For $\beta\in \mathbb{R}^{W}$, define the correlation for the two base densities: $h_{i}$ and $h_{j}$,
\begin{equation}\nonumber
\rho_{W}(i,j)=\frac{<h_{i},h_{j}>}{\|h_{i}\|\|h_{j}\|}, \quad i,j=1,\cdots,W.
\end{equation}

Our results will be established under the local coherence condition, we define the maximal local coherence as:
\begin{equation}\nonumber
    \rho(\beta) = \max_{i \in I(\beta)} \max_{j \neq i} |\rho_{W}(i,j)|
\end{equation}

It is easy to see that it measures the separation of the variables in the set $I(\beta)$ from one another and from the rest, and the degree of separation is measured in terms of the size of the correlation coefficients. However, the regular condition introduced by this coherence may be too strong, it may exclude the cases where the ``correlation'' can be relatively large for a small number of pairs $(i, j)$ and almost zero for otherwise, for instance. So we consider the definition of cumulative local coherence given by \cite{Bunea10}:
\begin{equation}\nonumber
\rho_{*}(\beta)=\sum_{i\in I(\beta)}\sum_{j>i}|\rho_{W}(i,j)|.
\end{equation}
Define
$$H(\beta)=\max_{j\in I(\beta)}\frac{\omega_{j}}{v(\delta/2)\|h_{j}\|},~~F=\max_{1\leq j\leq W}\frac{v(\delta/2)\|h_{j}\|}{\tilde{\omega}_{j}}=\max_{1\leq j\leq W}\frac{\|h_{j}\|}{2\sqrt{2}L_{j}},$$
where ${\tilde{\omega}_{j}}:=2\sqrt{2}L_{j}v(\delta/2)$.

By using the definition of $\rho_{*}(\beta)$ and the notations above, we present the key result of this paper which lays the foundation for the oracle inequality of the estimated mixed coefficients.

\begin{theorem}\label{theorem 1}
Under the assumption (H.1)-(H.5) and $c=\frac{\min_{1\leq j\leq W}\{\tilde{\omega}_{j}\}}{B}$, then the true base functions $\{h_{j}\}_{j=1}^{W}$ satisfies cumulative local coherence assumption
\begin{equation}\label{eq:5}
12FH(\beta)\rho_{*}(\beta)\sqrt{W(\beta)}\leq \gamma,
\end{equation}
and all $0<\gamma\leq1$, we have the following oracle inequality,
\begin{align*}
&~~~~\|h_{\hat{\beta}}-h\|^{2}+\frac{\alpha_{opt1}(1-\gamma)}{(\alpha_{opt1}-1)}\sum_{j=1}^{W}\tilde{\omega}_{j}|\hat{\beta}_{j}-\beta_{j}|+\frac{\alpha_{opt1}}{\alpha_{opt1}-1}\sum_{j=1}^{W}c(\hat{\beta}_{j}-\beta_{j})^{2}\\
&\le \frac{\alpha_{opt1}+1}{\alpha_{opt1}-1}\|h_{\beta}-h\|^{2}+\frac{18\alpha_{opt1}^{2}}{\alpha_{opt1}-1}H^{2}(\beta)v^{2}(\delta/2)W(\beta)
\end{align*}
 with probability $1-\delta$ at least, where $\alpha_{opt1}=1+\sqrt{1+\frac{\|h_{\beta}-h\|^{2}}{9H^{2}(\beta)v^{2}(\delta/2)W(\beta)}}$.
\end{theorem}
It is worthy to note that here we use $\sqrt{W(\beta)}$ instead of $W(\beta)$, and the later is used in \cite{Bunea10}. The reason of the phenomenon is quite clean actually: we introduce the $\ell_2$ penalty and derive our result in a more general frame while \cite{Bunea10} derive their result under restriction $\gamma = 1$. Now, let us address the sparse Gram matrix ${\rm{\psi}}_{W}=(<h_{i},h_{j}>)_{1\leq i,j\leq W}$ with a small number of non-zero element in off-diagonal positions, define $\psi_{W}(i,j)$ as the element $(i,j)$-th of position $\psi_{W}$. The condition \eqref{eq:5} in Theorem {\ref{theorem 1}} can be transformed to the condition $$12SH(\beta)\sqrt{W(\beta)}\leq \gamma,$$
where the number $S$ is called the sparse index of matrix ${\rm{\psi}}_{W}$ which is defined as follows:
\begin{equation}\nonumber
S=|\{(i,j):i,j\in\{1,\cdots,W\},i>j~\textrm{and}~\psi_{W}(i,j)\neq 0\}|,
\end{equation}
where $|A|$ is the number of elements of set $A$.

Sometimes the assumption \eqref{eq:5} does not necessarily imply the positive definiteness of ${\rm{\psi}}_{W}$.  Next we give similar oracle inequality that is valid under the hypothesis that the Gram matrix ${\rm{\psi}}_{W}$ is positive definite.
\begin{theorem}\label{theorem 2}
Under the assumption (H.1)-(H.5) and Gram matrix ${\rm{\psi}}_{W}$ are positive definite with minimum eigenvalues greater than or equal to $\lambda_{W}>0$. For all $\beta\in \mathbb{R}^{W}$, we have the following oracle inequality with  probability at least $1-\delta$,
\begin{align*}
&~~~~\|h_{\hat{\beta}}-h\|^{2}+\frac{\alpha_{opt2}}{\alpha_{opt2}-1}\sum_{j=1}^{W}\tilde{\omega}_{j}|\hat{\beta}_{j}-\beta_{j}|+\frac{\alpha_{opt2}}{\alpha_{opt2}-1}\sum_{j=1}^{W}c(\hat{\beta}_{j}-\beta_{j})^{2}\\
&\leq \frac{\alpha_{opt2}+1}{\alpha_{opt2}-1}\|h_{\beta}-h\|^{2}+\frac{576\alpha_{opt2}^{2}}{\alpha_{opt2}-1}\frac{{{G}}}{{{\lambda _W}}}v^{2}(\delta/2),
\end{align*}
where $G=G(\beta):=\sum_{j\in I(\beta)}L_{j}^{2}$ and $\alpha_{opt2}=1+\sqrt{1+\frac{\|h_{\beta}-h\|^{2}}{288\frac{{{G}}}{{{\lambda _W}}}v^{2}(\delta/2)}}$.
\end{theorem}
\textbf{Remark:} The argument and result of Theorem {\ref{theorem 1}} in this paper is more refined than the conclusion of Theorem {\ref{theorem 1}} in \cite{Bunea10} for Lasso by letting $\gamma=1/2$ and $c=0$. In addition, Theorem {\ref{theorem 1}} and Theorem {\ref{theorem 2}} of this paper respectively give the optimal $\alpha$ value of the density estimation oracle inequalities, namely $\alpha_{opt1}$, $\alpha_{opt2}$. It provides potentially sharper bound for the $\ell_1$-estimation error bound.

Next, we will present $\ell_1$-estimation error for the estimator $\hat\beta$ by (\ref{eq:9}), and the weights are defined by \eqref{eq:3.1}.

For technical reason, we consider that $\|h_{j}\|=1$ for all $j$ in \eqref{eq:9}, i.e. the based functions are normalized. This normalization mimics the covariates' standardization procedure when we do some penalized estimation in generalized linear models. For simplicity, we put $L:=\max_{1\leq j\leq W}L_{j}$.

For any other choice of $v(\delta/2)$ greater than or equal to $ \sqrt {\frac{1}{n}\log \frac{2W}{\delta }}$, the conclusions of Section {\ref{section 3}} are valid with a high probability. It imposes restriction on the predictive performance of CSDE. As pointed out in \cite{Bunea08}, for the $\ell_{1}$-penalty in the regression, the adjustment sequence $\omega_{j}$ required for the correct selection is usually larger than the adjustment sequence $\omega_{j}$ that produces a good prediction. The selection of the mixed density shown below is also true. Specifically, we will take the value $\beta = \beta^*$ and $v=v(\delta/2W)=\sqrt{\frac{\log(2W^{2}/\delta)}{n}}$ then $\alpha_{opt1},\alpha_{opt2}=2$, in below we give the corollary of Theorem {\ref{theorem 1}},{\ref{theorem 2}}.

\begin{corollary}\label{corollary 1}
Given the same conditions as Theorem {\ref{theorem 1}} with $\|h_{j}\|=1$ for all $j$ , let $\alpha_{opt1}=2$, we have the following $\ell_1$-estimation error oracle inequality:
\begin{align}\label{eq:4.3}
\sum_{j=1}^{W}|\hat{\beta}_{j}-\beta_{j}^{*}|\leq\frac{72\sqrt{2}v(\delta/2W)W(\beta^{*})}{1-\gamma}\frac{(L+L_{\min})^{2}}{L_{\min}}
\end{align}
with probability $1-\delta/W$ at least, where $L_{\min}=\min_{1\leq j\leq W}L_{j}$.
\end{corollary}
\begin{corollary}\label{corollary 2}
Given the same conditions as Theorem {\ref{theorem 2}} with $\|h_{j}\|=1$ for all $j$, let $\alpha_{opt2}=2$, we have the following $\ell_1$-estimation error oracle inequality, with probability $1-\delta/W$ at least,
\begin{align*}
\sum_{j=1}^{W}|\hat{\beta}_{j}-\beta_{j}^{*}|\leq\frac{288\sqrt{2}v(\delta/2W)G^{*}}{L_{\min}\lambda_{W}},
\end{align*}
where $G^{*}=\sum_{j\in I_*}L_{j}^{2}$.
\end{corollary}

If the number $W(\beta^{*})$ of the mixed indicator elements is much smaller than $\sqrt{n}$, then the inequality (\ref{eq:4.3}) guarantees that the estimated $\hat{\beta}$ is close to the true $\beta^{*}$, and the $\ell_1$-estimation error will be presented in the numerical simulation in Section {\ref{section 4}}. Our results of Corollary \ref{corollary 1} and {\ref{corollary 2}} are non-asymptotic for any $W$ and $n$.  The oracle inequalities are guiders for us to find an optimal tuning parameter with order $O(\sqrt {\frac{{\log W}}{n}} )$ for sharper estimation error and better prediction performance. This is also an intermediate and crucial result which leads to the main result of correctly identifying the mixture components in Section {\ref{section 3.2}}. In the following section, we turn to cope with the identification of $I_*$. The selection of correct components is derived by the proposed oracle inequality for the weighted $\ell_1$+ $\ell_2$ penalty.

\subsection{Correct support identification of mixture models}\label{section 3.2}

In this section, we will study results on the support recovery of our CSDE estimator. There are a few versions of support recovery while most of the results are the consistency of $\ell_1$-error and prediction error. Here we borrow the framework of \cite{Bunea08} and \cite{Zhang20}, they give many proof techniques to deal with the correct support identification in linear models by $\ell_{1}+\ell_{2}$ regularization. Let $\hat{I}$ be the set of indicators consisting of non-zero elements of $\hat{\beta}$ in the given (\ref{eq:9}). In other words, $\hat{I}$ is an estimate of the true support set $I(\beta^*):=I_{*}$. We will study that for a given $0<\varepsilon<1$, $$P(\hat{I}=I(\beta^*))\geq 1-\varepsilon$$
under some mild conditions.

In order to identify the $I_{*}$ consistently, we need more assumptions about some special correlation conditions than $\ell_1$-error consistency.

\textbf{Condition (A)}:
\[{\rho _*}({\beta ^*}) \le \frac{{L{L_{\min }}{\lambda _W}}}{{{\rm{288}}{G^*}}}.\]
Moreover, we need an additional condition that the minimal signal should be higher than a threshold level and quantified by the order of tuning parameter. We state it as follows:

\textbf{Condition (B)}:
\begin{align*}
\min_{j\in I^{*}}|\beta_{j}^{*}|\geq4\sqrt{2}v(\frac{\delta}{2W})L,
\end{align*}
where $v(\frac{\delta}{2W}): = \sqrt {\frac{1}{n}\log \frac{2W^2}{\delta }}$.

When performing simulation, condition (B) is the theoretical guarantee that the smallest magnitude of $\beta_j$ must be greater than a threshold value as a minimal signal condition. It is also called Beta-min condition, see \cite{Buhlmann2011}.

\begin{theorem}\label{theorem 3}
Let $0<\delta<\frac{1}{2}$ be a given value and define $\epsilon_{k}:=|E[h_{k}(X_{1})]-E[h_{k}(Z_{1})]|$. Assume that both condition A and B are true, and give the same conditions as Corollary {\ref{corollary 2}}, then
\begin{align*}
P(\hat{I}=I_{*})\geq 1-\left(4W\left(\frac{\delta}{2W^{2}}\right)^{(1-\epsilon_{k}^{*})^2}+2\delta\right),
\end{align*}
where $\epsilon_{k}^{*}=\epsilon_{k}/\sqrt{2}v(\delta/2W)L$.
\end{theorem}

Under the beta-min condition, the support estimation is very close to the true the support of $\beta_{j}^{*}$. The probability of the event $\{\hat{I}=I_{*}\}$ is high when $W$ is growing. The $\hat{\beta}$ recovers the correct support with the probability at least $ 1-\left(4W\left(\frac{\delta}{2W^{2}}\right)^{(1-\epsilon_{k}^{*})^2}+2\delta\right)$. The result is non-asymptotic, it is true for any fixed $W$ and $n$. Similar conclusion about support consistency, see Theorem 6 of {\cite{Zhang20}}.

\section{Simulation and Real Data Analysis}\label{section 4}

\cite{Bunea10} suppose the spades estimation to deal with the samples for sparse mixture density, and they also derive an algorithm to complement their theoretical result. Their findings successfully handle the high-dimensional adaptive density estimation in some degree. However, their algorithm is costly and unstable. In this section, we deal with the tuning parameter directly and compare our CSDE method with the SPADES method (AdaLasso) in \cite{Bunea10} and other similar methods. In all cases here, we fix $n = 100$ for $W = 81, 131, 211, 321$, which is known as the dimension of the unknown parameter $\beta^*$. The performance of each estimator is evaluated by $\ell_1$-estimation error and total variation (TV) distance between the estimator and the true value of $\beta^*$. The total variation (TV) error is defined as:
$${\rm{TV}}( h_{\beta^*}, h_{\hat{\beta}} ) = \int | h_{\beta^*}(x) - h_{\hat{\beta}}(x)| dx.$$

\subsection{Tuning parameter selection}

In \cite{Bunea10}, the $\lambda_1$ is chosen by coordinate descent method, while the mixture weights are detected by general bisection method (GBM). But in our article, the optimal weights can be computed directly. So it is much easier to carry out than \cite{Bunea10}. The $\ell_1$-penalty term $ \sum_{j = 1} ^ W \omega_j |\beta_j|$ with optimal weights is defined by
\begin{equation*}
\omega_{k}:=2\sqrt{2}L_{k}v(\delta/2)+cB,
\end{equation*}
where $L_j = \| h_j \|_{\infty}$, which usually can be computed easily for a continuous $h_j$.

For a discrete base density $\{h_{j}\}_{j=1}^{W}$, it can be estimated as the following approximation by using concentration inequalities from Exercise 4.3.3 of \cite{chow1997probability}:
\begin{equation}
|\text{med}(X) - E(X)| \le \sqrt {2Var(X)},
\end{equation}
\begin{equation}
    \bar{x} \approx x_{med} \left( 1 + O ( n ^ {-1} ) \right) \approx h ^ {-1} (L_j) \left( 1 + O ( n ^ {-1} )  \right),
\end{equation}
where $\bar{x}$ and $x_{med}$ represent the sample mean and sample median respectively in each simulation, then we only need to select the $\lambda_1$ and $c = \lambda_2$, and they can be detected by nesting coordinate descent method. Besides, the precision level is assigned as $\xi = 0.001$ in our simulation.

\subsection{Multi-modal distributions}
First, we examine our method in a multi-modal Gaussian model that is similar to the first model in \cite{Bunea10}. However, our mixture Gaussian has a different variance, which leads the meaningful weights to our estimation. The model is assigned as follow:
\begin{equation}
    h_{\beta}^* (x) = \sum_{j = 1} ^ W \beta_j^*\phi \left( x| a\cdot j, \sigma _j \right).
\end{equation}
%In practice, the difficulty is to detect the sparse high-dimensional parameter in the sample, which often is not sufficiently large in size. The multi-modal will also increase the complication of the weak signals.
We choose $a = 0.5, n=100$ and:
\begin{equation}
    \beta^* = \left( \mathbf{0}_8^{T}, 0.2, \mathbf{0}_{10}^{T}, 0.1, \mathbf{0}_5^{T}, 0.1, \mathbf{0}_{10}^{T}, 0.1, \mathbf{0}_{10}^{T}, 0.1, \mathbf{0}_5^{T}, 0.15, \mathbf{0}_{10}^{T}, 0.15, \mathbf{0}_{10}^{T}, 0.1, \mathbf{0}_{W - 76}^{T} \right)^{T}.
\end{equation}
%The variances of the Gaussian distribution are also the signal:
with
\begin{equation}
    \sigma = \left( \mathbf{1}_{20}^{T}, \mathbf{0.8}_6^{T}, \mathbf{0.6}_{11}^{T}, \mathbf{0.4}_{11}^{T}, \mathbf{0.6}_{6}^{T}, \mathbf{0.8}_{11}^{T}, \mathbf{1.2}_{W - 76}^{T}\right) ^{T}.
\end{equation}

An acceptable measurement error $e_j (x)$ satisfied $E e_j(X) = 0$ is chosen as:
\begin{equation}
    h_j(x) + e_j(x) \, \sim \, N \left(a\cdot j, 1.1 \sigma_j^2 \right).
\end{equation}
Then we use the sample $x_1, \ldots, x_n \ \text{i.i.d.} \, \sum_{j = 1} ^ W \beta_j^{*} [h_j(x) + e_j(x)]$ to estimate $\beta^{*}$. We replicate the simulation $N = 100$ times. The results of simulation are presented in Table \ref{table1}, we can see our method has the more and more excellent performances as the $W$ increases which matches the non-asymptotical results in the previous section. The best performance is far better than the other three method when $W = 321$. It's worthy to note that the better approximation following the increasing of $W$, matching the \label{eq:kkt1} and Theorem {\ref{theorem 3}} in our previous section.
\begin{table}[!htb]\large
    \centering
\emph{
\caption{\emph{The mean and standard deviation of the errors in the four estimators of $\beta^*$ under $N = 100$ simulations, with $n = 100$. The quasi-optimal $\lambda_2$ is $c = 0.002$ for Enet, while $c = 0.027$  is for the adaptive method.\\}}\label{table1}
\begin{tabular}{ccccc}
\hline \hline
         & $W$                   & $\lambda_1$                & $L_1$ error        & $TV$ error       \\ \hline
Lasso    & \multirow{4}{*}{81}  & \multirow{2}{*}{0.065} & 2.133 (2.467)   & 1.137 (1.115)  \\ \cline{4-5}
Enet     &                      &                        & 2.061 (1.439)   & 1.114 (0.805)  \\ \cline{3-5}
AdaLasso &                      & \multirow{2}{*}{0.053} & 1.922 (2.211)   & 1.258 (1.296)  \\ \cline{4-5}
CSDE  &                      &                        & 2.191 (4.812)   & 1.405 (2.329)  \\ \hline
Lasso    & \multirow{4}{*}{131} & \multirow{2}{*}{0.068} & 2.032 (0.985)   & 1.352 (0.712)  \\ \cline{4-5}
Enet     &                      &                        & 2.236 (2.498)   & 1.409 (1.056)  \\ \cline{3-5}
AdaLasso &                      & \multirow{2}{*}{0.056} & 1.880 (2.644)   & 0.972 (1.204)  \\ \cline{4-5}
CSDE  &                      &                        & 1.635 (0.342)   & 0.863 (0.402)  \\ \hline
Lasso    & \multirow{4}{*}{211} & \multirow{2}{*}{0.071} & 2.572 (4.187)   & 1.605 (2.702)  \\ \cline{4-5}
Enet     &                      &                        & 2.061 (1.883)   & 1.353 (1.516)  \\ \cline{3-5}
AdaLasso &                      & \multirow{2}{*}{0.058} & 1.764 (1.041)   & 0.832 (0.610)  \\ \cline{4-5}
CSDE  &                      &                        & 1.648 (0.168)   & 0.791 (0.415)  \\ \hline
Lasso    & \multirow{4}{*}{321} & \multirow{2}{*}{0.074} & 2.120 (2.842)   & 1.146 (1.115)  \\ \cline{4-5}
Enet     &                      &                        & 10.173 (82.753) & 7.839 (67.887) \\ \cline{3-5}
AdaLasso &                      & \multirow{2}{*}{0.061} & 2.106 (4.816)   & 0.818 (1.565)  \\ \cline{4-5}
CSDE  &                      &                        & 1.623 (0.085)   & 0.634 (0.199)  \\ \hline \hline
\end{tabular}
}
\end{table}

\begin{figure}[H]
    \centering
    \includegraphics[width=1.1\textwidth]{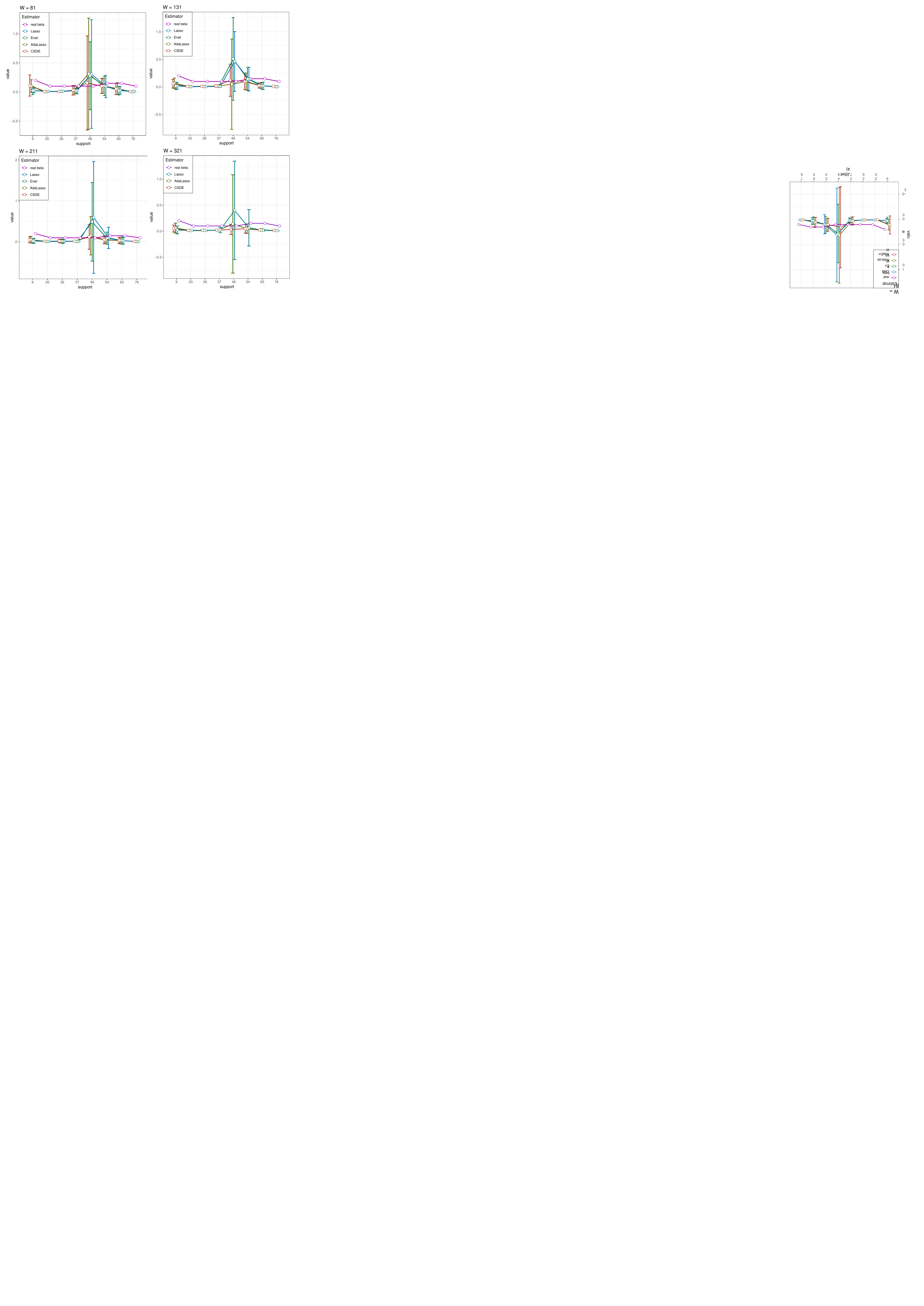} % Figure image
    \caption{\emph{The estimated support of $\beta^*$ by the four types of estimators, and the $W$ is varying. The circles represent the means of the estimators under the four specific approaches, while the half of the vertical lines mean the standard deviations. }} % Figure caption
    \label{figure1}
\end{figure}

\begin{figure}[!htp]
    \centering
    \includegraphics[width=1.05\textwidth]{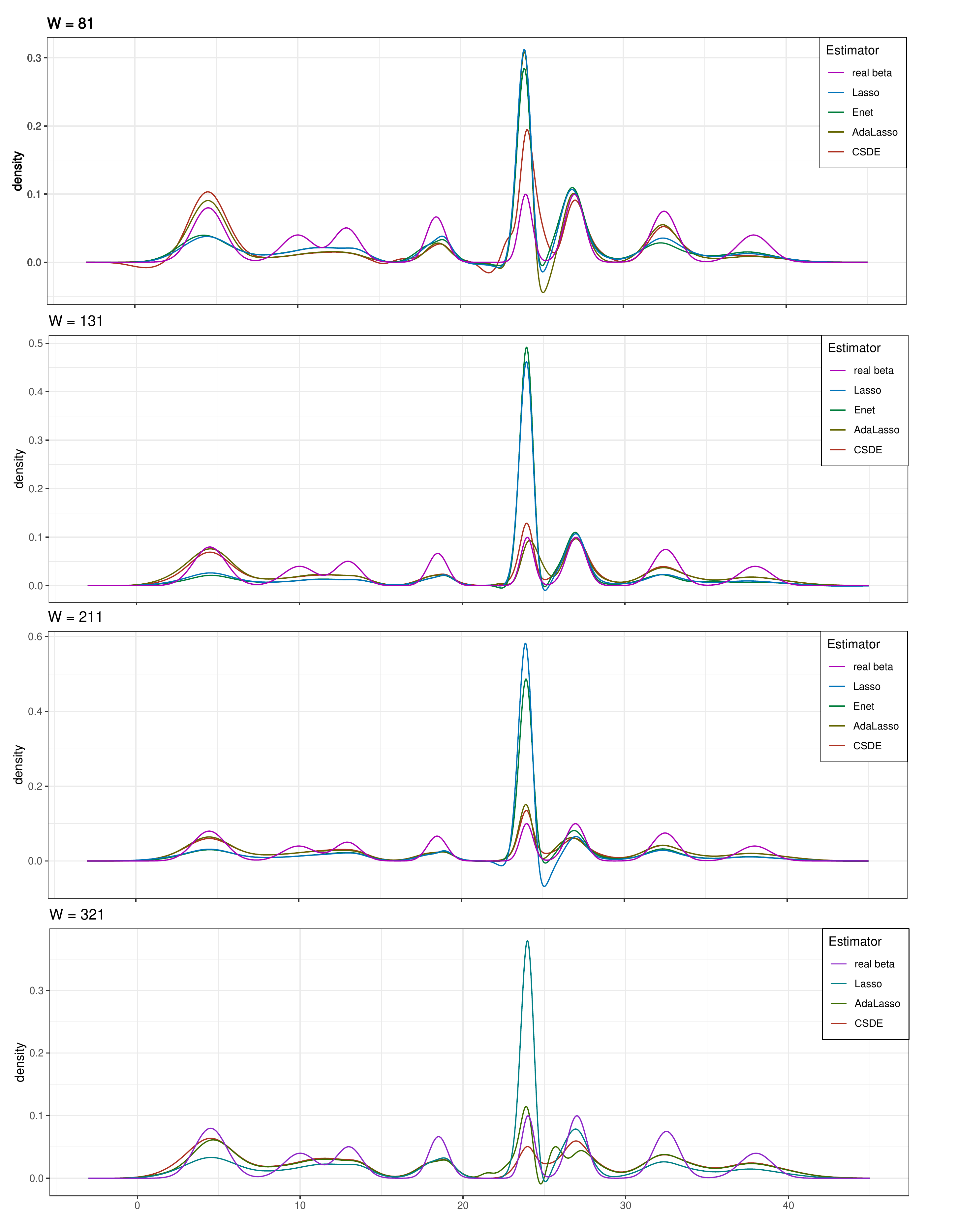} % Figure image
    \caption{\emph{The density map of the four estimators' means and the true denses. The result of $Enet$ in $W = 321$  is not be shown due to its exactly poor performance.}} % Figure caption
    \label{figure2}
\end{figure}

We plot the solution path to compare the performance of the four estimators in $\beta_j \in I(\beta)$ for every $W$ in Figure \ref{figure1} (the result of Enet in $W = 321$ is not be shown due to its poor performance.). These figures also provide strong support for the above analysis. Meanwhile, we plot the probability densities of the several estimators and the true density to complement the visual sensory of the advantage in our method in Figure \ref{figure2}, in which the powerful competency of detecting the multi-mode is shown (whereas other methods only find the most strong signal, ignoring other meaningful but relatively weak signals).

\subsection{Mixture of Poisson distributions}

In the second set of our simulations, we study the mixture of discrete distribution: the mixture Poisson distribution
\begin{equation}
    h_{\beta^*} (x) = \sum_{j = 1} ^ W \beta_j^* p \left( x| \lambda_j = a\cdot j \right),
\end{equation}
where $a = 0.1$, and:
\begin{equation}
    \beta^* = \left( \mathbf{0}_8^{T}, 0.2, \mathbf{0}_{10}^{T}, 0.1, \mathbf{0}_5^{T}, 0.1, \mathbf{0}_{10}^{T}, 0.1, \mathbf{0}_{10}^{T}, 0.1, \mathbf{0}_5^{T}, 0.15, \mathbf{0}_{10}^{T}, 0.15, \mathbf{0}_{10}^{T}, 0.1, \mathbf{0}_{W - 75}^{T} \right)^{T}.
\end{equation}

The adjusted weights are calculated by \eqref{eq:9}, and in discrete distributions, we define $\langle f,g\rangle  = \sum\limits_{k = 1}^\infty  {f(k)g(k)}$. Meanwhile, the Poisson distribution with measurement errors can be treated as a negative binomial distribution. Namely:
\begin{equation}
    p \left( x| \lambda_i \right) + e_i(x) \ \sim NB \left(r, \frac{\lambda_i}{\lambda_i + r} \right),
\end{equation}
where $\lambda_i$ is the mean for $i$-th observation and $r$ is the common dispersion parameter. Here, $i=1,2,\cdots,n=100$.

For a practical $r$, we choose $r = 6$, which leads an increment $\lambda_ij ^2 / r$ in the variance. Similarly, we replicate each simulation to estimate the parameter $N = 100$ times with the sample coming from the mixture negative binomial distribution above. The result is shown in Table \ref{table2}. The result is actually akin to that in the previous mixture Gaussian distribution, while the strong performance of our method is shown clearly when $W$ is considerable.

\begin{table}[H]\large
    \centering
\emph{
    \caption{\emph{The mean and standard deviation of the errors in the four estimators of $\beta$ under $N = 100$ simulations. The $\lambda_2$ is chosen as $c = 0.005$ for Enet, while $c = 0.203$ for the adaptive method.\\}}
    \label{table2}
\begin{tabular}{ccccc}
\hline \hline
         & $W$                   & $\lambda_1$                & $L_1$ error        & $TV$ error       \\ \hline
Lasso    & \multirow{4}{*}{81}  & \multirow{2}{*}{0.048} & 1.796 (0.006)   & 0.002 (0.001)  \\ \cline{4-5}
Enet     &                      &                        & 1.796 (0.006)   & 0.002 (0.001)  \\ \cline{3-5}
AdaLasso &                      & \multirow{2}{*}{0.138} & 1.811 (0.013)   & 0.002 (0.005)  \\ \cline{4-5}
CSDE  &                      &                        & 1.806 (0.008)   & 0.003 (0.005)  \\ \hline
Lasso    & \multirow{4}{*}{131} & \multirow{2}{*}{0.051} & 1.828 (0.006)   & 0.003 (0.001)  \\ \cline{4-5}
Enet     &                      &                        & 1.830 (0.009)   & 0.004 (0.002)  \\ \cline{3-5}
AdaLasso &                      & \multirow{2}{*}{0.145} & 1.880 (0.006)   & 0.002 (0.005)  \\ \cline{4-5}
CSDE  &                      &                        & 1.854 (0.006)   & 0.002 (0.004)  \\ \hline
Lasso    & \multirow{4}{*}{211} & \multirow{2}{*}{0.053} & 1.935 (0.010)   & 0.005 (0.003)  \\ \cline{4-5}
Enet     &                      &                        & 2.061 (0.014)   & 0.007 (0.008)  \\ \cline{3-5}
AdaLasso &                      & \multirow{2}{*}{0.152} & 1.935 (0.008)   & 0.005 (0.003)  \\ \cline{4-5}
CSDE  &                      &                        & 1.861 (0.005)   & 0.003 (0.002)  \\ \hline
Lasso    & \multirow{4}{*}{321} & \multirow{2}{*}{0.055} & 1.927 (0.031)   & 0.005 (0.002)  \\ \cline{4-5}
Enet     &                      &                        & 2.123 (0.026)   & 0.009 (0.009) \\ \cline{3-5}
AdaLasso &                      & \multirow{2}{*}{0.158} & 1.938 (0.008)   & 0.005 (0.003)  \\ \cline{4-5}
CSDE  &                      &                        & 1.852 (0.002)   & 0.002 (0.001)  \\ \hline \hline
\end{tabular}
}
\end{table}

\subsection{Low dimensional mixture model}
Surprisingly, our method has more competitive efficacy than some popular methods (such as EM algorithm), even the dimension $W$ is relatively small. To see this, we introduce the following numerical experiments to estimate the weights of the low dimensional mixed Gaussian model: the samples $x_{1},\cdots,x_{n}$  come from the model:
\begin{equation}\nonumber
h_{\beta^*}(x)=\sum_{j=1}^{W} \beta_{j}^* \phi\left(x | \mu_{j}, \sigma_{j}\right).
\end{equation}
The updating equation for EM algorithm in $t$-th step is:
\begin{equation}\nonumber
\omega_{i j}^{(t)}=\frac{p_{j}^{(t)} \phi\left(x_{i} ; \mu_{t}, \sigma_{t}\right)}{\sum_{s=1}^{W} p_{s}^{(t)} \phi\left(x_{i} ; \mu_{s}, \sigma_{s}\right)}, \quad \beta_{j}^{(t+1)}=\frac{\sum_{i=1}^{W} \omega_{i j}^{(t)}}{\sum_{i=1}^{n} \sum_{j=1}^{W} \omega_{i j}^{(t)}}.
\end{equation}
Here we consider two scenarios:
\begin{align*}
&{\text { (1) } W=6, \beta=(0.3,0,0,0.3,0,0.4)^{T}, \mu=(0,10,20,30,40,50)^{T}, \sigma=(1,2,3,4,5,6)^{T}} \\
&\text { (2) } W=7, \beta=(0.1,0,0,0.8,0,0,0.1)^{T}, \mu=(0,1,2,3,4,5,6)^{T}, \sigma=(0.3,0.2,0.2,0.1,0.2,0.2,0.3)^{T}.
\end{align*}

For each scenario $n=50$, and the fitter levels (cessation level) in the EM approach and our method are both $\xi=10^{-4}$. A well-advised initial value in the EM approach is the equal weight.

We replicate the simulation $N = 100$ times, and the optimal tuning parameters stem from the CV (so under each simulation they are not the same, albeit they are very close to each other). The result can be seen in Table \ref{table3}.

\begin{table}[!htb]
    \centering

\caption{\emph{The eventual simulations result.\\}}
\label{table3}

\begin{tabular}{ccccc}
\hline \hline
                            &         & $L_1$ error             & $TV$ error                \\ \hline
\multirow{2}{*}{Scenario 1} & EM      & 0.2547758(0.1217747)    & 0.2045831(0.09807962)     \\ \cline{3-4}
                            & CSDE & 0.2055284(0.1449865)    & 0.1852503(0.104226)       \\ \cline{1-4}
\multirow{2}{*}{Scenario 2} & EM      & 0.1109372(0.05491049)   & 0.1107578(0.05491172)     \\ \cline{3-4}
                            & CSDE & 0.1090387(0.03715101)   & 0.108296(0.03683571)      \\ \hline
\hline
\end{tabular}

\end{table}

\subsection{Real data examples}

\begin{figure}[H]
    \centering
    \includegraphics[width=1\textwidth]{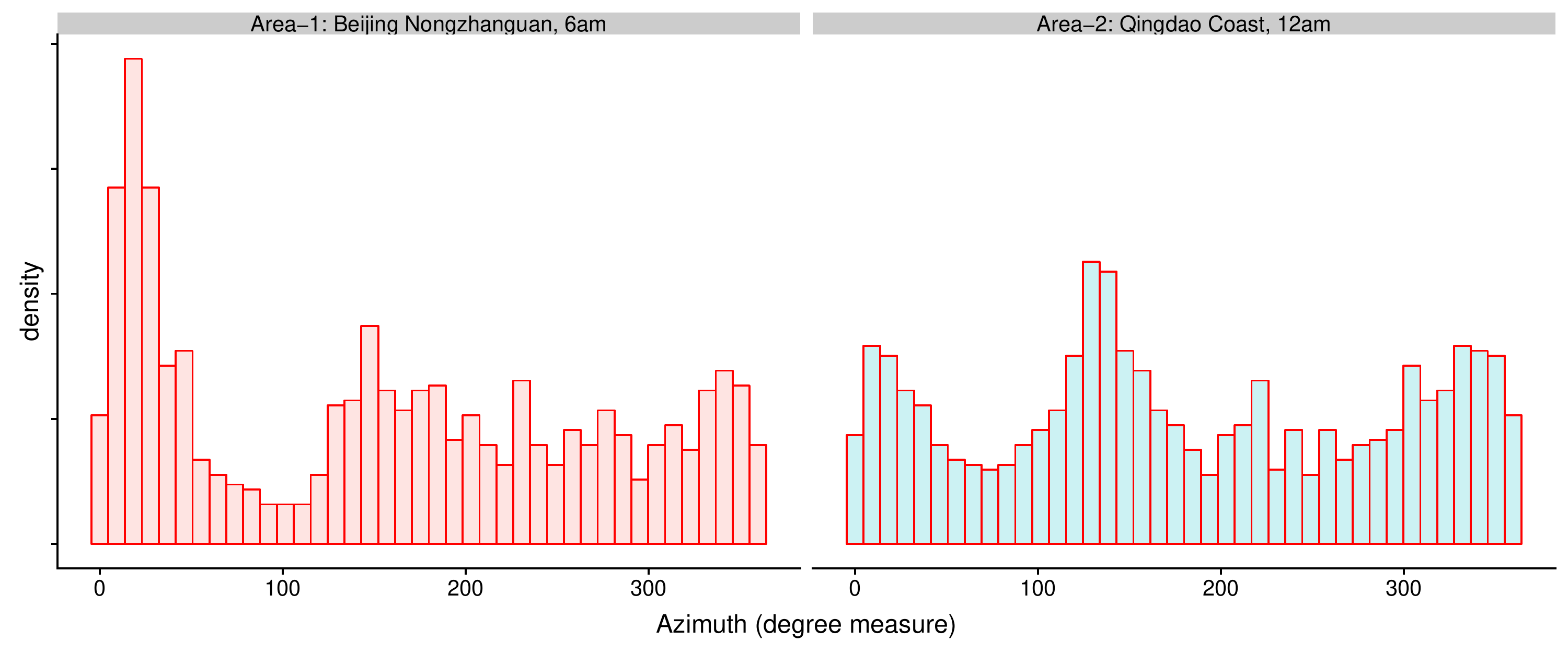} % Figure image
    \caption{\emph{The sample histogram of the azimuth in Beijing Nongzhanguan at 6am and Qingdao Coast at 12am.}} % Figure caption
    \label{figure3}
\end{figure}

Practically, we consider using our method to estimate some densities in the environmental science field. In the area of meteorology, wind, which is mercurial, has been an advisable object to study for a long time. Take notice of the speed of the wind at one specific location maybe not diverse, so we will use the wind's azimuth angle with a more sparse density at two locations in China. Concerning many types of research about the estimated density for wind existed, so there is a possibility to use our approach to cope with some difficulties in meteorology science.

There have been some very credible meteorological data sets. We would use the ERA5 hourly data in \cite{Hersbach2018} to continue our analysis. We would like to choose a continental area and a coastal area in China which refers to Beijing Nongzhanguan and Qingdao Coast, respectively. The location of these two area are: $[116.3125 \degree E, 116.4375 \degree E] \times [39.8125 \degree N , 39.8125 \degree N]$. Take notice that the wind in one day may be highly correlated, therefore, using the data at a specific time point of each day in a consecutive period as i.i.d. samples are reasonable. The sample histogram of 6 am in Beijing Nongzhanguan and at 12 am Qingdao Coast is shown in Figure \ref{figure3}. Here we use the data from 2013/01/01 to 2015/12/12.

\begin{figure}[!htb]
    \centering
    \includegraphics[width=1.1\textwidth]{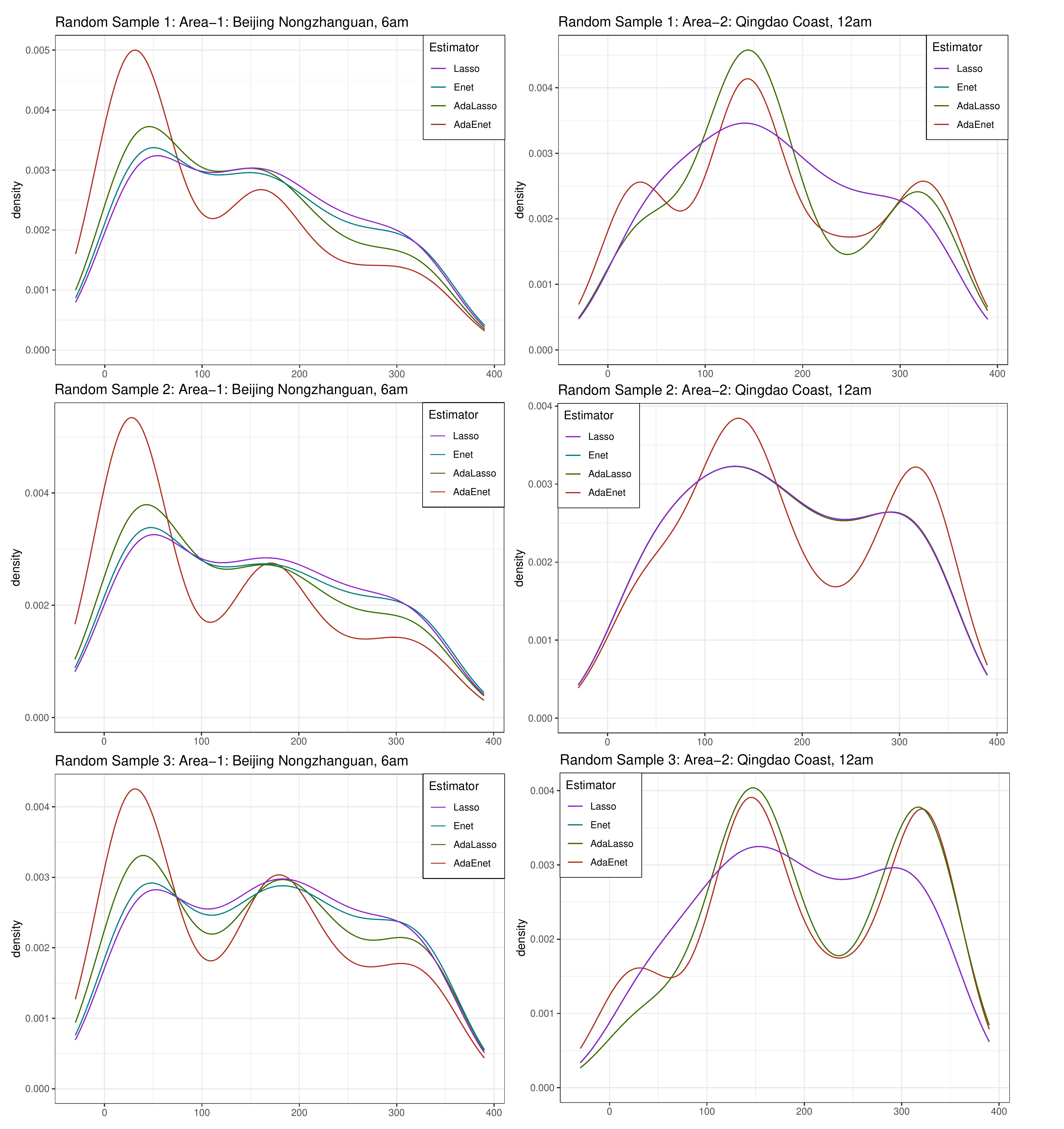} % Figure image
    \caption{\emph{The density map of the four estimators' approaches for the three random subsamples.}} % Figure caption
    \label{figure4}
\end{figure}

As we can see, their density does be multi-peak (we used 1095 samples). Now we can use our approach to estimate the multi-mode densities based on a relatively small size of samples, which is only a tiny part of the whole data from 2013/01/01 to 2015/12/12. Because one year has nearly 360 days, we may assume that every day is a latent factor that forms the base density. So the model is designed as
\begin{equation}
    h_{\beta^*} (x) = \sum_{j = 1} ^ {360} \beta_j^* \phi \left( x| \mu_j, \sigma_j \right)
\end{equation}
with the mean and variance parameters:
\begin{equation}
    \mu = (1, 2, \ldots, 360)^{T}, \quad \sigma = t \cdot \mathbf{1}_{360}^{T},
\end{equation}
where the $t$ is seen as the bandwidth (or tuning parameter). With the different sub-samples, the values are computed are different.

Another critical issue is how to choose the tuning parameter $\lambda_1$ and $\lambda_2$. Then we apply the cross-validation criterion, namely choosing $\lambda_i$ minimizing the difference of two estimators derived from the separated samples in random dichotomy.

Now start to construct the samples for estimating procedure. Assume that an observatory wants to figure some information about the wind in the two areas. However, it doesn't have intact data due to the limited budget at its inception. The only samples it has are $6 \sim 8$ days' information each month for the 2 areas and these days scatter randomly. This imperfect data does increase the challenge of estimating the trustworthy density, and we would compare our method with other previous methods, in which appraising the difference between the complete data sample histogram and the estimated density under each method is for the evaluation. Notice that the samples are only a tiny part of the whole data, so the $n$ is relatively small. The small sample and large dimension setting coincide with the non-asymptotical theory provided in the previous section. The estimating density has been shown in Figure \ref{figure4}.

Evidently, in this practical application, we can see that our method vindicates its more efficient estimating performance and stability from its propinquity of the complete sample histogram, namely the efficacious capacity of detecting the shape of multi-mode density, and the stronger inclination to bear resemblance to each other subsample (although some subtle nuances do exist by reason of the different subsample).

\section{Summary and Discussions}\label{section 5}
The paper deals with the deconvolution problem using Lasso type methods: the observations are $X_{1}, \cdots, X_{n}$ which are independent and generated from $X_{i}=Z_{i}+\varepsilon_{i}$; and the goal is to estimate the unknown density $h$ of the $Z_{i}$. we assume that the function $h$ can be written as
$h(\cdot)=h_{\beta^{*}}(\cdot)=\sum_{j=1}^{W} \beta_{j}^{*} h_{j}(\cdot)$ based on some functions $\{h_{j}\}_{j=1}^{W}$ from a specific dictionary and propose to estimate the coefficients of this decomposition with the Elastic-net method. For this estimator, we show that under some classical assumptions on the model such as coherence of the Gram matrix, finite sample bounds both for the estimation and the prediction errors valid with a relatively high probability can be obtained. Moreover, under a beta-min condition, we prove a variable selection consistency result. An extensive numerical study is also conducted. The following estimation problem is also similar to the CSDE.

%In this paper, we study the finite sample properties of sparse mixed density estimation models by deriving non-asymptotic oracle inequalities of the weighted $\ell_1$+$\ell_2$  penalized estimator $\hat\beta$. Our purpose is to estimate the vector $\beta^{*}$ by adapting the weighted $\ell_1$+$\ell_2$  penalty to this unknown sparsity of $\beta^{*}$ and then to identify $I_{*}$ with a high probability.

\emph{Aggregate density estimator with measurement errors.}  Based on the idea of model average, our aim is to aggregate some candidate density models $h_{1},\cdots,h_{W}$ based on the data $\{X_{i}\}_{i=1}^{n}$ containing measurement errors. It means we need to construct a new aggregated estimator as the convex combination of $h_{1},\cdots,h_{W}$, which is approximately the best among $h_{1},\cdots,h_{W}$. The aggregation we consider here in the form of $h_{\hat{\beta}}$ by appropriately chosen the weight vector $(\beta_{1},\cdots,\beta_{W})\in \mathbb{R}^{W}$.

For the future study, it is also interesting and meaningful to do hypothesis testing about the coefficients $\beta^{*} \in \mathbb{R}^{W}$ in sparse mixture models. For a general function $h: \mathbb{R}^{W} \mapsto \mathbb{R}^{m}$ and a nonempty closed set $\Omega \in \mathbb{R}^{m}$, we can consider
$$
H_{0}: h\left(\beta^{*}\right) \in \Omega \text { vs. } H_{1}: h\left(\beta^{*}\right) \notin \Omega.
$$

\appendix
 \section{Technical Details}
For convenience, we first give a preliminary lemma and proof, now define the random variables
\begin{equation}\nonumber
M_{j}=\frac{1}{n}\sum_{i=1}^{n}\left\{h_{j}(X_{i})-Eh_{j}(X_{i})\right\}.
\end{equation}
Consider the event $\cal{E}$ by
\begin{equation}\nonumber
{\cal{E}}=\bigcap_{j=1}^{W}\{2|M_{j}|\leq \tilde{\omega}_{j}\},
\end{equation}
where $\tilde\omega_{k}:=2\sqrt{2}L_{k}\sqrt{\frac{1}{n}\log \frac{W}{\delta/2}}=:2\sqrt{2}L_{k}v(\delta/2)$.

Then we have the following lemma which is cornerstone for the proofs in below.
\begin{lemma}\label{lemma 3}
Suppose that $\mathop {\max }\limits_{1 \le j \le W} {L_j} < \infty $ and $c=\frac{\min_{1\leq j\leq W}\{\tilde{\omega}_{j}\}}{B}$, for any $\beta\in \mathbb{R}^{W}$, on the event $\mathcal{E}$, we have
\begin{align}
\|h_{\hat{\beta}}-h\|^{2}+\sum_{j=1}^{W}{\tilde{\omega}}_{j}|\hat{\beta}_{j}-\beta_{j}|+\sum_{j=1}^{W}c(\hat{\beta}_{j}-\beta_{j})^{2}
\leq \|h_{\beta}-h\|^{2}+6\sum_{j\in I(\beta)}{\omega}_{j}|\hat{\beta}_{j}-\beta_{j}|.
\end{align}
\end{lemma}
\subsection*{Proof of Lemma {\ref{lemma 3}}}
According to the definition of $\hat{\beta}$, for any $\beta\in \mathbb{R}^{W}$, we find
\begin{equation}
-\frac{2}{n}\sum_{i=1}^{n}h_{\hat{\beta}}(X_{i})+\|h_{\hat{\beta}}\|^{2}+2\sum_{j=1}^{W}\omega_{j}|\hat{\beta}_{j}|+c\sum_{j=1}^{W}\hat{\beta}_{j}^{2}\leq -\frac{2}{n}\sum_{i=1}^{n}h_{\beta}(X_{i})+\|h_{\beta}\|^{2}+2\sum_{j=1}^{W}\omega_{j}|\beta_{j}|+c\sum_{j=1}^{W}\beta_{j}^{2}.
\end{equation}
Then
\begin{align}\label{lm:then}
\|h_{\hat{\beta}}\|^{2}-\|h_{\beta}\|^{2}\leq
\frac{2}{n}\sum_{i=1}^{n}h_{\hat{\beta}}(X_{i})-\frac{2}{n}\sum_{i=1}^{n}h_{\beta}(X_{i})+2\sum_{j=1}^{W}\omega_{j}|\beta_{j}|
-2\sum_{j=1}^{W}\omega_{j}|\hat{\beta}_{j}|+c\sum_{j=1}^{W}\beta_{j}^{2}-c\sum_{j=1}^{W}\hat{\beta}_{j}^{2}.
\end{align}
Note that
\begin{align}\label{lm:note}
\|h_{\hat{\beta}}-h\|^{2}&=\|h_{\hat{\beta}}-h_{\beta}+h_{\beta}-h\|^{2}\nonumber\\
&=\|h_{\hat{\beta}}-h_{\beta}\|^{2}+\|h_{\beta}-h\|^{2}+2<h_{\beta}-h,h_{\hat{\beta}}-h_{\beta}>\nonumber\\
&=\|h_{\beta}-h\|^{2}-2<h,h_{\hat{\beta}}-h_{\beta}>+2<h_{\beta},h_{\hat{\beta}}-h_{\beta}>+\|h_{\hat{\beta}}-h_{\beta}\|^{2}\nonumber\\
&=\|h_{\beta}-h\|^{2}-2<h,h_{\hat{\beta}}-h_{\beta}>+\|h_{\hat{\beta}}\|^{2}-\|h_{\beta}\|^{2}.
\end{align}
Combining \eqref{lm:then} and \eqref{lm:note}, we obtain
\begin{align}\label{lm:comb}
\|h_{\hat{\beta}}-h\|^{2}&\leq\|h_{\beta}-h\|^{2}+2\sum_{j=1}^{W}\omega_{j}|\beta_{j}|-2\sum_{j=1}^{W}\omega_{j}|\hat{\beta}_{j}|+c\sum_{j=1}^{W}\beta_{j}^{2}-c\sum_{j=1}^{W}\hat{\beta}_{j}^{2}\nonumber\\
&-2<h,h_{\hat{\beta}}-h_{\beta}>+\frac{2}{n}\sum_{i=1}^{n}h_{\hat{\beta}}(X_{i})-\frac{2}{n}\sum_{i=1}^{n}h_{\beta}(X_{i}).
\end{align}
According to the definition of $h_{\beta}(x)$, it gives $h_{\beta}(x)=\sum_{j=1}^{W}\beta_{j}h_{j}(x)$ with $\beta=(\beta_{1},\cdots,\beta_{W})$. For the 3 terms in \eqref{lm:comb}, we have
\begin{align*}
&~~~~-2<h,h_{\hat{\beta}}-h_{\beta}>+\frac{2}{n}\sum_{i=1}^{n}h_{\hat{\beta}}(X_{i})-\frac{2}{n}\sum_{i=1}^{n}h_{\beta}(X_{i})\\
&=2\cdot\frac{1}{n}\sum_{i=1}^{n}\left(\sum_{j=1}^{W}\hat{\beta}_j h_{j}(X_{i})-\sum_{j=1}^{W}\beta_{j}h_{j}(X_{i})\right)-2E(h_{\beta '}-h_{\beta})(X_{i})|_{\beta '=\hat{\beta}}\\
&=2\sum_{j=1}^{W}\frac{1}{n}\sum_{i=1}^{n}h_{j}(X_{i})(\hat{\beta}_{j}-\beta_{j})-2\sum_{j=1}^{W}E[h_{j}(X_{i})](\hat{\beta}_{j}-\beta_{j})\\
&=2\sum_{j=1}^{W}\left(\frac{1}{n}\sum_{i=1}^{n}h_{j}(X_{i})-E[h_{j}(X_{i})]\right)(\hat{\beta}_{j}-\beta_{j}).
\end{align*}
Then
\begin{align*}
\|h_{\hat{\beta}}-h\|^{2}&\le \|h_{\beta}-h\|^{2}+2\sum_{j=1}^{W}\left(\frac{1}{n}\sum_{i=1}^{n}h_{j}(X_{i})-E[h_{j}(X_{i})]\right)(\hat{\beta}_{j}-\beta_{j})\\
&+2\sum_{j=1}^{W}\omega_{j}|\beta_{j}|-2\sum_{j=1}^{W}\omega_{j}|\hat{\beta}_{j}|+c\sum_{j=1}^{W}\beta_{j}^{2}-c\sum_{j=1}^{W}\hat{\beta}_{j}^{2}.
\end{align*}
Conditioning on the event $\mathcal{E}$, we have
\begin{align*}
\|h_{\hat{\beta}}-h\|^{2}\leq \|h_{\beta}-h\|^{2}+\sum_{j=1}^{W}\tilde{\omega}_{j}|\hat{\beta}_{j}-\beta_{j}|+2\sum_{j=1}^{W}\omega_{j}(|\beta_{j}|-|\hat{\beta}_{j}|)+c\sum_{j=1}^{W}(\beta_{j}^{2}-\hat{\beta}_{j}^{2}).
\end{align*}
Adding $\sum_{j=1}^{W}\tilde{\omega}_{j}|\hat{\beta}_{j}-\beta_{j}|+c\sum_{j=1}^{W}(\beta_{j}-\hat{\beta}_{j})^{2}$ to both sides of the above inequality, it gives
\begin{align*}
&~~~~\|h_{\hat{\beta}}-h\|^{2}+\sum_{j=1}^{W}\tilde{\omega}_{j}|\hat{\beta}_{j}-\beta_{j}|+c\sum_{j=1}^{W}(\beta_{j}-\hat{\beta}_{j})^{2}\\
&\leq \|h_{\beta}-h\|^{2}+2\sum_{j=1}^{W}\tilde{\omega}_{j}|\hat{\beta}_{j}-\beta_{j}|+2\sum_{j=1}^{W}\omega_{j}(|\beta_{j}|-|\hat{\beta}_{j}|)+c\sum_{j=1}^{W}(\beta_{j}^{2}-\hat{\beta}_{j}^{2})+c\sum_{j=1}^{W}(\beta_{j}-\hat{\beta}_{j})^{2}.
\end{align*}
Note that
\begin{align*}
&~~~~c[\sum_{j=1}^{W}(\beta_{j}^{2}-\hat{\beta}_{j}^{2})+\sum_{j=1}^{W}(\beta_{j}-\hat{\beta}_{j})^{2}]=c[\sum_{j=1}^{W}(\beta_{j}^{2}-\hat{\beta}_{j}^{2}+\beta_{j}^{2}-2\beta_{j}\hat{\beta}_{j}+\hat{\beta}_{j}^{2})]\\
&=2c\sum_{j=1}^{W}\beta_{j}(\beta_{j}-\hat{\beta}_{j})=2c\sum_{j\in I(\beta)}\beta_{j}(\beta_{j}-\hat{\beta}_{j})\leq2cB\sum_{j\in I(\beta)}|\beta_{j}-\hat{\beta}_{j}|\leq2\sum_{j\in I(\beta)}\tilde{\omega}_{j}|\beta_{j}-\hat{\beta}_{j}|,
\end{align*}
where the last inequality is due to the assumption $c=\frac{\min_{1\leq j\leq W}\{\tilde{\omega}_{j}\}}{B}\leq\frac{\tilde{\omega}_{j}}{B}$.

So we obtain
\begin{align*}
&~~~~\|h_{\hat{\beta}}-h\|^{2}+\sum_{j=1}^{W}{\tilde{\omega}}_{j}|\hat{\beta}_{j}-\beta_{j}|+c\sum_{j=1}^{W}(\hat{\beta}_{j}-\beta_{j})^{2}\\
&\leq\|h_{\beta}-h\|^{2}+2\sum_{j=1}^{W}\tilde{\omega}_{j}|\hat{\beta}_{j}-\beta_{j}|+2\sum_{j=1}^{W}\omega_{j}(|\beta_{j}|-|\hat{\beta}_{j}|)+2\sum_{j\in I(\beta)}{\tilde{\omega}}_{j}|\hat{\beta}_{j}-\beta_{j}|\\
&\leq\|h_{\beta}-h\|^{2}+2\sum_{j=1}^{W}{\omega}_{j}|\hat{\beta}_{j}-\beta_{j}|+2\sum_{j=1}^{W}\omega_{j}(|\beta_{j}|-|\hat{\beta}_{j}|)+2\sum_{j\in I(\beta)}{\omega}_{j}|\hat{\beta}_{j}-\beta_{j}|,
\end{align*}
where the last inequality follows from $\tilde{\omega}_{j} \le {\omega}_{j}$ for all $j$.

We know that when $j\in I(\beta)$, $\beta_{j}\neq0$, when $j\notin I(\beta)$, $\beta_{j}=0$. Considering $|\beta_{j}|-|\hat{\beta}_{j}|\leq|\hat{\beta}_{j}-\beta_{j}|$ for all $j$, we have
\begin{align*}
2\sum_{j=1}^{W}&{\omega}_{j}|\hat{\beta}_{j}-\beta_{j}|+2\sum_{j=1}^{W}\omega_{j}(|\beta_{j}|-|\hat{\beta}_{j}|)\le 4\sum_{j\in I(\beta)}{\omega}_{j}|\hat{\beta}_{j}-\beta_{j}|.
\end{align*}
Then
\begin{align*}
\|h_{\hat{\beta}}-h\|^{2}+\sum_{j=1}^{W}{\tilde{\omega}}_{j}|\hat{\beta}_{j}-\beta_{j}|+c\sum_{j=1}^{W}(\hat{\beta}_{j}-\beta_{j})^{2}&\leq\|h_{\beta}-h\|^{2}+4\sum_{j\in I(\beta)}{\omega}_{j}|\hat{\beta}_{j}-\beta_{j}|+2\sum_{j\in I(\beta)}{\omega}_{j}|\hat{\beta}_{j}-\beta_{j}|\\
&=\|h_{\beta}-h\|^{2}+6\sum_{j\in I(\beta)}{\omega}_{j}|\hat{\beta}_{j}-\beta_{j}|.
\end{align*}
This completes the proof.

\subsection*{The Proof of Theorems}

According to the construction of $\tilde{\omega}_{j}=2\sqrt{2}L_{j}\sqrt{\frac{1}{n}\log\frac{2W}{\delta}}$ in the formula (\ref{eq:3.1}), the sum of the independent random variables $\zeta_{ij}=h_{j}(X_{i})-Eh_{j}(X_{i})$ is determined by Hoeffding's inequality, and $|\zeta_{ij}|\leq2L_{j}$. We obtain
\begin{align*}
P(\mathcal{E}^{c})=P\left(\bigcup_{j=1}^{W}\{2|M_{j}|>\tilde{\omega}_{j}\}\right)&\leq\sum_{j=1}^{W}P(2|M_{j}|>\tilde{\omega}_{j})\leq2\sum_{j=1}^{W}\exp\left(-\frac{2n^{2}\cdot\tilde{\omega}_{j}^{2}/4}{4nL_{j}^{2}}\right)\\
&=2\sum_{j=1}^{W}\exp\left(-\log\frac{2W}{\delta}\right)=2W\cdot\frac{\delta}{2W}=\delta.
\end{align*}

\subsection*{Proof of Theorem {\ref{theorem 1}}}
By using Lemma {\ref{lemma 3}}, we need an upper bound on $\sum_{j\in I(\beta)}\omega_{j}|\hat{\beta}_{j}-\beta_{j}|$. For easy notation, let
\begin{equation}\nonumber
q_{j}=\hat{\beta}_{j}-\beta_{j}, \quad Q(\beta)=\sum_{j\in I(\beta)}|q_{j}|\|h_{j}\|, \quad Q=\sum_{j=1}^{W}|q_{j}|\|h_{j}\|.
\end{equation}
According to the definition of $H(\beta)$, that is, $H(\beta)=\max_{j\in I(\beta)}\frac{{\omega}_{j}}{v(\delta/2)\|h_{j}\|}$, we have
\begin{equation}\label{eq:wv}
\sum_{j\in I(\beta)}{\omega}_{j}|\hat{\beta}_{j}-\beta_{j}|\leq v(\delta/2)H(\beta)Q(\beta).
\end{equation}
Let $Q_{*}(\beta):=\sqrt{\sum_{j\in I(\beta)}q_{j}^{2}\|h_{j}\|^{2}}$. Using the definition of $h_{\beta}(x)$, we obtain
\begin{align*}
Q_{*}^{2}(\beta)=\sum_{j\in I(\beta)}q_{j}^{2}\|h_{j}\|^{2}&=\|h_{\hat{\beta}}-h_{\beta}\|^{2}-\sum_{i,j\notin I(\beta)}q_{i}q_{j}<h_{i},h_{j}>\\
&-(2\sum_{i\notin I(\beta)}\sum_{j\in I(\beta)}q_{i}q_{j}<h_{i},h_{j}>
+\underset{i,j\in I(\beta), i\neq j}{\sum\sum}q_{i}q_{j}<h_{i},h_{j}>).
\end{align*}
As $i,j\notin I(\beta)$, $\beta_{i}=\beta_{j}=0$, it is easy to see,
\begin{equation}\nonumber
\underset{i,j\notin I(\beta)}{\sum\sum}<h_{i},h_{j}>q_{i}q_{j}\geq0.
\end{equation}
Observe that
\begin{align*}
&~~~~2\sum_{i\notin I(\beta)}\sum_{j\in I(\beta)}q_{i}q_{j}<h_{i},h_{j}>+\underset{i,j\in I(\beta), i\neq j}{\sum\sum}q_{i}q_{j}<h_{i},h_{j}>\\
&=2\sum_{i\notin I(\beta)}\sum_{j\in I(\beta)}q_{i}q_{j}<h_{i},h_{j}>+2\underset{i,j\in I(\beta), j>i}{\sum\sum}q_{i}q_{j}<h_{i},h_{j}>=2\underset{i\in I(\beta), j>i}{\sum\sum}q_{i}q_{j}<h_{i},h_{j}>.
\end{align*}
By the definitions of $\rho_{W}(i,j)$ and $\rho_{*}(\beta)$, then
\begin{align*}
Q_{*}^{2}(\beta)&\leq\|h_{\hat{\beta}}-h_{\beta}\|^{2}+2\underset{i\in I(\beta), j>i}{\sum\sum}|q_{i}||q_{j}|\|h_{i}\|\|h_{j}\|\frac{<h_{i},h_{j}>}{\|h_{i}\|\|h_{j}\|}\\
&\leq\|h_{\hat{\beta}}-h_{\beta}\|^{2}+2\rho_{*}(\beta)\max_{i\in I(\beta),j>i}|q_{i}|\|h_{i}\||q_{j}|\|h_{j}\|.
\end{align*}
In fact,
\begin{align*}
\max_{i\in I(\beta)}{|q_{i}|\|h_{i}\|}\leq\sqrt{\sum_{j\in I(\beta)}q_{j}^{2}\|h_{j}\|^{2}}=Q_{*}(\beta),~~\max_{i\in I(\beta),j>i}|q_{j}|\|h_{j}\|\leq\sum_{j=1}^{W}|q_{j}|\|h_{j}\|.
\end{align*}
Thus
\begin{align}\label{eq:6}
Q_{*}^{2}(\beta)&\leq\|h_{\hat{\beta}}-h_{\beta}\|^{2}+2\rho_{*}(\beta)Q_{*}(\beta)\sum_{j=1}^{W}|q_{j}|\|h_{j}\|=\|h_{\hat{\beta}}-h_{\beta}\|^{2}+2\rho_{*}(\beta)Q_{*}(\beta)Q.
\end{align}
By (\ref{eq:6}), we can get
\begin{align*}
Q_{*}^{2}(\beta)-2\rho_{*}(\beta)Q_{*}(\beta)Q-\|h_{\hat{\beta}}-h_{\beta}\|^{2}\leq0.
\end{align*}

In order to find the upper bound of $Q_{*}(\beta)$, apply the properties of the quadratic inequality to the above formula,
\begin{align}\label{eq:7}
Q_{*}(\beta)\leq\rho_{*}(\beta)Q+\sqrt{\rho^{2}_{*}(\beta)Q^{2}+\|h_{\hat{\beta}}-h_{\beta}\|^{2}}&\leq \rho_{*}(\beta)Q+[\rho_{*}(\beta)Q+\|h_{\hat{\beta}}-h_{\beta}\|]\nonumber\\
\leq2\rho_{*}(\beta)Q+\|h_{\hat{\beta}}-h_{\beta}\|.
\end{align}
Note that $W(\beta)=|I(\beta)|=\sum_{j=1}^{W}I(\beta_{j}\neq0)$, employing Cauchy-Schwarz inequalities, we have
\[W(\beta )\sum\limits_{j \in I(\beta )} | {q_j}|^2\|{h_j}\|{^2} = \sum\limits_{j \in I(\beta )} {I^2}(j \in I(\beta )) \sum\limits_{j \in I(\beta )} | {q_j}{|^2}\|{h_j}\|{^2} \ge {(\sum\limits_{j \in I(\beta )} {{I({\rm{\{ }}j \in I(\beta )})} |} {q_j}|{\|{h_j}\|})^2 = {Q^2}(\beta ).\]
Then
\begin{align*}
Q_{*}^{2}(\beta)=\sum_{j\in I(\beta)}|q_{j}|^{2}\|h_{j}\|^{2}\geq Q^{2}(\beta)/W(\beta).
\end{align*}
Combined with (\ref{eq:7}), we can get $
Q(\beta)/\sqrt{W(\beta)}\leq Q_{*}(\beta)\leq2\rho_{*}(\beta)Q+\|h_{\hat{\beta}}-h_{\beta}\|$. Therefore,
\begin{align}\label{eq:8}
Q(\beta)\leq2\rho_{*}(\beta)\sqrt{W(\beta)}Q+\sqrt{W(\beta)}\|h_{\hat{\beta}}-h_{\beta}\|.
\end{align}

By Lemma {\ref{lemma 3}}, we have the following inequality established by probability exceeding $1-\delta$.
\begin{align*}
&~~~~\|h_{\hat{\beta}}-h\|^{2}+\sum_{j=1}^{W}{\tilde{\omega}}_{j}|\hat{\beta}_{j}-\beta_{j}|+\sum_{j=1}^{W}c(\hat{\beta}_{j}-\beta_{j})^{2}\\
&\leq \|h_{\beta}-h\|^{2}+6\sum_{j\in I(\beta)}{\omega}_{j}|\hat{\beta}_{j}-\beta_{j}|~~(\text{by}~\eqref{eq:wv})\\
&\leq \|h_{\beta}-h\|^{2}+6v(\delta/2)H(\beta)Q(\beta)~~(\text{by}~\eqref{eq:8})\\
&\leq\|h_{\beta}-h\|^{2}+6v(\delta/2)H(\beta)[2\rho_{*}(\beta)\sqrt{W(\beta)}\sum_{j=1}^{W}|q_{j}|\|h_{j}\|+\sqrt{W(\beta)}\|h_{\hat{\beta}}-h_{\beta}\|]\\
&= {\left\| {{h_\beta } - h} \right\|^2} + 12v(\delta/2)H(\beta )\rho_{*}(\beta )\sqrt{W(\beta )}\sum\limits_{j = 1}^W {{\tilde{\omega} _j}|{{\hat \beta }_j} - {\beta _j}|\frac{{\left\| {{h_j}} \right\|}}{{{\tilde{\omega} _j}}}} + 6v(\delta/2)H(\beta )\sqrt {W(\beta )}\| {{h_{\hat \beta }} - {h_\beta }} \|\\
&\leq  {\left\| {{h_\beta } - h} \right\|^2} + 12FH(\beta )\rho_{*}(\beta )\sqrt{W(\beta )}\sum\limits_{j = 1}^W {{\tilde{\omega} _j}|{{\hat \beta }_j} - {\beta _j}|}  + 6v(\delta/2)H(\beta )\sqrt {W(\beta )} \| {{h_{\hat \beta }} - {h_\beta }} \|\\
&\leq  {\left\| {{h_\beta } - h} \right\|^2} + \gamma\sum\limits_{j = 1}^W {{\tilde{\omega} _j}|{{\hat \beta }_j} - {\beta _j}|}  + 6v(\delta/2)H(\beta )\sqrt {W(\beta )}\| {{h_{\hat \beta }} - {h_\beta }} \|,
\end{align*}
where the second last inequality follows from the definition of $F:=\mathop {\max }\limits_{1 \le j \le W} \frac{{v(\delta /2)\left\| {{h_j}} \right\|}}{{{{\tilde \omega }_j}}}$, and the last inequality is derived by the assumption $12FH(\beta )\rho_{*}(\beta )\sqrt{W(\beta )} \le \gamma ,(0 < \gamma  \le 1).$

Further, we can find that, with probability at least $1-\delta$,
\begin{align*}
&~~~~\|h_{\hat{\beta}}-h\|^{2}+(1-\gamma)\sum_{j=1}^{W}{\tilde{\omega}}_{j}|\hat{\beta}_{j}-\beta_{j}|+\sum_{j=1}^{W}c(\hat{\beta}_{j}-\beta_{j})^{2}\\
&\leq \|h_{\beta}-h\|^{2}+6v(\delta/2)H(\beta)\sqrt{W(\beta)}\|h_{\hat{\beta}}-h_{\beta}\|\\
&=\|h_{\beta}-h\|^{2}+6v(\delta/2)H(\beta)\sqrt{W(\beta)}\|h_{\hat{\beta}}-h+h-h_{\beta}\|\\
&\leq\|h_{\beta}-h\|^{2}+6v(\delta/2)H(\beta)\sqrt{W(\beta)}\|h_{\hat{\beta}}-h\|+6v(\delta/2)H(\beta)\sqrt{W(\beta)}\|h-h_{\beta}\|.
\end{align*}
Using the elementrary inequality $2st\leq s^{2}/\alpha+\alpha t^{2}~(s,t\in \mathbb{R},\alpha>1)$ to the last two terms of the above inequality, it yields
\begin{align*}
&2\{3v(\delta/2)H(\beta)\sqrt{W(\beta)}\}\|h_{\hat{\beta}}-h\|\leq\alpha\cdot9v^{2}(\delta/2)H^{2}(\beta)W(\beta)+\|h_{\hat{\beta}}-h\|^{2}/\alpha,\\
&2\{3v(\delta/2)H(\beta)\sqrt{W(\beta)}\}\|h_{\beta}-h\|\leq\alpha\cdot9v^{2}(\delta/2)H^{2}(\beta)W(\beta)+\|h_{\beta}-h\|^{2}/\alpha.
\end{align*}
Thus
\begin{align*}
&~~~~\|h_{\hat{\beta}}-h\|^{2}+(1-\gamma)\sum_{j=1}^{W}{\tilde{\omega}}_{j}|\hat{\beta}_{j}-\beta_{j}|+\sum_{j=1}^{W}c(\hat{\beta}_{j}-\beta_{j})^{2}\\
&\leq\|h_{\beta}-h\|^{2}+18\alpha v^{2}(\delta/2)H^{2}(\beta)W(\beta)+\|h_{\hat{\beta}}-h\|^{2}/\alpha+\|h_{\beta}-h\|^{2}/\alpha.
\end{align*}
Simplifying, we have
\begin{align}\label{eq:oracle}
&~~~~\|h_{\hat{\beta}}-h\|^{2}+\frac{\alpha(1-\gamma)}{(\alpha-1)}\sum_{j=1}^{W}\tilde{\omega}_{j}|\hat{\beta}_{j}-\beta_{j}|+\frac{\alpha}{\alpha-1}\sum_{j=1}^{W}c(\hat{\beta}_{j}-\beta_{j})^{2}\nonumber\\
&\leq \frac{\alpha+1}{\alpha-1}\|h_{\beta}-h\|^{2}+\frac{18\alpha^{2}}{\alpha-1}H^{2}(\beta)v^{2}(\delta/2)W(\beta),~\alpha > 1,0 < \gamma  \le 1.
\end{align}
Optimizing $\alpha$ to obtain the sharp upper bounds for the above oracle inequality
\[{\alpha _{opt1}}: = \mathop {\arg \min }\limits_{\alpha  > 1} \left\{ {\frac{\alpha+1}{\alpha-1}\|h_{\beta}-h\|^{2}+\frac{18\alpha^{2}}{\alpha-1}H^{2}(\beta)v^{2}(\delta/2)W(\beta)} \right\}=1+\sqrt{1+\frac{\|h_{\beta}-h\|^{2}}{9H^{2}(\beta)v^{2}(\delta/2)W(\beta)}}\]
by the first order condition.

So far, Theorem {\ref{theorem 1}} is proved by substituting ${\alpha _{opt1}}$ into \eqref{eq:oracle}.\\

\subsection*{Proof of Theorem {\ref{theorem 2}}}
By the minimal eigenvalue assumption for $\psi_{W}$, we have
\begin{align}\label{eq:leasteigen}
\|h_{\beta}\|^{2}=\|\sum_{j=1}^{W}\beta_{j}h_{j}(x)\|^{2}=\beta^{T}\psi_{W}\beta\geq\lambda_{W}\|\beta\|^{2}\geq\lambda_{W}\sum_{j\in I(\beta)}\beta_{j}^{2}.
\end{align}
Using the definition of ${\omega}_{j}$ and assumption ${L_{\min }}: = \mathop {\min }\limits_{1 \le j \le W} {L_j} > 0$, thus
\[{\omega _j} = 2{L_j}\left( {\sqrt {\frac{{2\log (2W/\delta )}}{n}} {\rm{ + }}\frac{{cB}}{{{2L_j}}}} \right) \le 2{L_j}\left( {\sqrt {\frac{{2\log (2W/\delta )}}{n}} {\rm{ + }}\frac{{cB}}{{{2L_{\min }}}}} \right).\]
Since $cB=\tilde{\omega}_{\min}=2\sqrt{2}L_{\min}v(\delta/2)$ and $v(\delta/2)=\sqrt{\frac{\log(2W/\delta)}{n}}$, we have
\[{\omega _j}\leq4\sqrt{2}L_{j}v(\delta/2).\]

Let $G(\beta)=\sum_{j\in I(\beta)}L_{j}^{2}$, we can get by the Cauchy-Schwartz inequality
\begin{align}\label{eq:leasteigenh}
6\sum_{j\in I(\beta)}{\omega}_{j}|\hat{\beta}_{j}-\beta_{j}|&\leq24\sqrt{2}v(\delta/2)\sum_{j\in I(\beta)}L_{j}|\hat{\beta}_{j}-\beta_{j}|\nonumber\\
&\leq24\sqrt{2}v(\delta/2)\sqrt{\sum_{j\in I(\beta)}L_{j}^{2}}\sqrt{\sum_{j\in I(\beta)}(\hat{\beta}_{j}-\beta_{j})^{2}}\leq 24\sqrt{2}v(\delta/2)\sqrt {\frac{{{G}(\beta)}}{{{\lambda _W}}}}\|h_{\hat{\beta}}-h_{\beta}\|,
\end{align}
where the last inequality above is due to
$$\|h_{\hat{\beta}}-h_{\beta}\|^{2}=\underset{1\leq i,j\le W}{\sum\sum}(\hat{\beta}_{i}-\beta_{i})(\hat{\beta}_{j}-\beta_{j})<h_{i},h_{j}> \ge\lambda_{W}\sum_{j\in I(\beta)}(\hat{\beta}_{j}-\beta_{j})^{2}$$
from \eqref{eq:leasteigen}.

Let $b(\beta):=12\sqrt{2}v(\delta/2)\sqrt {\frac{{{G}(\beta)}}{{{\lambda _W}}}} $,  Lemma {\ref{lemma 2}} implies
\begin{align*}
\|h_{\hat{\beta}}-h\|^{2}+\sum_{j=1}^{W}{\tilde{\omega}}_{j}|\hat{\beta}_{j}-\beta_{j}|+\sum_{j=1}^{W}c(\hat{\beta}_{j}-\beta_{j})^{2}
&\leq \|h_{\beta}-h\|^{2}+2b(\beta)\|h_{\hat{\beta}}-h_{\beta}\|\\
&= \|h_{\beta}-h\|^{2}+2b(\beta)(\|h_{\hat{\beta}}-h+h-h_{\beta}\|)\\
&\leq \|h_{\beta}-h\|^{2}+2b(\beta)\|h_{\hat{\beta}}-h\|+2b(\beta)\|h_{\beta}-h\|.
\end{align*}

Using the inequality $2st\leq s^{2}/\alpha+\alpha t^{2} \quad (s,t\in R, \alpha>1)$ for the last two terms on the right side of the above inequality, we find
\begin{align*}
2b(\beta)\|h_{\hat{\beta}}-h\|+2b(\beta)\|h_{\beta}-h\|&\leq\|h_{\hat{\beta}}-h\|^{2}/\alpha+b^{2}(\beta)\alpha+\|h_{\beta}-h\|^{2}/\alpha+b^{2}(\beta)\alpha\\
&=\|h_{\hat{\beta}}-h\|^{2}/\alpha+\|h_{\beta}-h\|^{2}/\alpha+2b^{2}(\beta)\alpha.
\end{align*}
Thus
\begin{align*}
\|h_{\hat{\beta}}-h\|^{2}+\sum_{j=1}^{W}{\tilde{\omega}}_{j}|\hat{\beta}_{j}-\beta_{j}|+\sum_{j=1}^{W}c(\hat{\beta}_{j}-\beta_{j})^{2}
\leq\|h_{\beta}-h\|^{2}+\|h_{\hat{\beta}}-h\|^{2}/\alpha+\|h_{\beta}-h\|^{2}/\alpha+2b^{2}(\beta)\alpha.
\end{align*}
We have
\begin{align*}
\frac{\alpha-1}{\alpha}\|h_{\hat{\beta}}-h\|^{2}+\sum_{j=1}^{W}{\tilde{\omega}}_{j}|\hat{\beta}_{j}-\beta_{j}|+\sum_{j=1}^{W}c(\hat{\beta}_{j}-\beta_{j})^{2}
\leq\frac{\alpha+1}{\alpha}\|h_{\beta}-h\|^{2}+2\alpha b^{2}(\beta).
\end{align*}
Therefore
\begin{align*}
\|h_{\hat{\beta}}-h\|^{2}+\frac{\alpha}{\alpha-1}\sum_{j=1}^{W}{\tilde{\omega}}_{j}|\hat{\beta}_{j}-\beta_{j}|+\frac{\alpha}{\alpha-1}\sum_{j=1}^{W}c(\hat{\beta}_{j}-\beta_{j})^{2}&\leq\frac{\alpha+1}{\alpha-1}\|h_{\beta}-h\|^{2}+\frac{2\alpha^{2}}{\alpha-1} b^{2}(\beta)\\
&=\frac{\alpha+1}{\alpha-1}\|h_{\beta}-h\|^{2}+\frac{576\alpha^{2}}{\alpha-1}\frac{{{G}(\beta)}}{{{\lambda _W}}}v^{2}(\delta/2).
\end{align*}
To get the sharp upper bounds for the above oracle inequality, we optimize $\alpha$
\begin{align*}
~~~{\alpha _{opt2}}: &= \mathop {\arg \min }\limits_{\alpha  > 1} \left\{ \frac{\alpha+1}{\alpha-1}\|h_{\beta}-h\|^{2}+\frac{576\alpha^{2}}{\alpha-1}\frac{{{G}(\beta)}}{{{\lambda _W}}}v^{2}(\delta/2) \right\}=1+\sqrt{1+\frac{\|h_{\beta}-h\|^{2}}{288\frac{{{G}(\beta)}}{{{\lambda _W}}}v^{2}(\delta/2)}},
\end{align*}
by the first order condition.  This completes the proof of Theorem {\ref{theorem 2}}.\\

\subsection*{Proof of Corollory {\ref{corollary 1}}}
 Let $\tilde{\omega}_{\min}:=\min_{1\leq j\leq W}\tilde{\omega}_{j}$. We replace $v(\delta/2)$ in Theorem {\ref{theorem 1}} by the larger value $v(\delta/2W)$. Substitute $\beta=\beta^{*}$ in Theorem {\ref{theorem 1}}, we have
\begin{align*}
\frac{\alpha_{opt1}(1-\gamma)}{\alpha_{opt1}-1}\sum_{j=1}^{W}\tilde{\omega}_{j}|\hat{\beta}_{j}-\beta_{j}^{*}|\leq\frac{18\alpha_{opt1}^{2}}{\alpha_{opt1}-1}
H^{2}(\beta^{*})v^{2}(\delta/2W)W(\beta^{*})
\end{align*}
by $h=h_{\beta^{*}}$.
Since $\tilde{\omega}_{j}\geq\tilde{\omega}_{\min}$ for all $j$, we get
\begin{align*}
\sum_{j=1}^{W}|\hat{\beta}_{j}-\beta_{j}^{*}|\leq\frac{18\alpha_{opt1}}{1-\gamma}\cdot\frac{1}{\tilde{\omega}_{\min}}\cdot\max_{j\in I(\beta)}\frac{\omega_{j}^{2}}{\|h_{j}\|^{2}}\cdot W(\beta^{*}).
\end{align*}
In this case, $\alpha_{opt1}=2$ and $\|h_{j}\|=1$, thus
\begin{align*}
~~~~\|\hat{\beta}-\beta^{*}\|&\leq\frac{36}{1-\gamma}\cdot\max_{j\in I(\beta)}\frac{\omega_{j}^{2}}{\tilde{\omega}_{\min}}\cdot W(\beta^{*})\\
&=\frac{72\sqrt{2}v(\delta/2W)W(\beta^{*})}{1-\gamma}\max_{j\in I(\beta)}\frac{(L_{j}+L_{\min})^{2}}{L_{\min}}\leq\frac{72\sqrt{2}v(\delta/2W)W(\beta^{*})}{1-\gamma}\frac{(L+L_{\min})^{2}}{L_{\min}}
\end{align*}
from
\begin{center}
 $\tilde{\omega}_{\min}=2\sqrt{2}v(\delta/2W)L_{\min}$, $\omega_{j}^{2}=[2\sqrt{2}v(\delta/2W)]^{2}\left[L_{j}+\frac{\tilde{\omega}_{\min}}{2\sqrt{2}v(\delta/2W)}\right]^{2}=[2\sqrt{2}v(\delta/2W)]^{2}[L_{j}+L_{\min}]^{2}$.\\
\end{center}
This completes the proof of corollary {\ref{corollary 1}}.
\subsection*{Proof of Corollary {\ref{corollary 2}}}
Let $\beta=\beta^{*}$ in Theorem {\ref{theorem 2}}, with $\alpha_{opt2}=2$, we replace $v(\delta/2)$ in Theorem {\ref{theorem 2}} by the larger value $v(\delta/2W)$, then
\begin{align*}
\sum_{j=1}^{W}\tilde{\omega}_{\min}|\hat{\beta}_{j}-\beta_{j}^{*}|\leq\sum_{j=1}^{W}\tilde{\omega}_{j}|\hat{\beta}_{j}-\beta_{j}^{*}|
\leq\frac{576\alpha_{opt2}G^{*}}{\lambda_{W}}v^{2}(\delta/2W).
\end{align*}
By the definition of $\tilde{\omega}_{\min}$, we can obtain
\[\sum\limits_{j = 1}^W | {\hat \beta _j} - \beta _j^*| \le \frac{{576{\alpha _{opt2}}{G^*}{v^2}(\delta /2)}}{{{{\tilde \omega }_{\min }}{\lambda _W}}}{\rm{ = }}\frac{{576 \cdot {\rm{2}}{G^*}{v^2}(\delta /2W)}}{{2\sqrt 2 v(\delta /2W){L_{\min }}{\lambda _W}}}=\frac{{{\rm{288}}\sqrt 2 {G^*}v(\delta /2W)}}{{{L_{\min }}{\lambda _W}}}.\]
This concludes the proof of corollary {\ref{corollary 2}}.
\subsection*{Proof of Theorem {\ref{theorem 3}}}
The following lemma is  by virtue of KKT conditions. It derives a bound of $P(I_{*}\nsubseteqq\hat{I})$ which is easily analysed.
\begin{lemma}[Proposition 3.3 in \cite{Bunea08}]\label{lemma 4}
$$P(I_{*}\nsubseteqq\hat{I}) \leq W(\beta^{*})\max_{k\in I_{*}}P(\hat{\beta}_{k}=0~\textrm{and}~\beta_{k}^{*}\neq0).$$
\end{lemma}
To present the proof of Theorem {\ref{theorem 3}}, first we notice,
\begin{align*}
P(\hat{I}\neq I_{*})\leq P(I_{*}\nsubseteqq\hat{I})+P(\hat{I}\nsubseteqq I_{*}).
\end{align*}
Next, we control the probability on the right side of the above inequality.\\

For the control of $P(I_{*}\nsubseteqq\hat{I})$, by Lemma {\ref{lemma 4}}, it remains to find $P(\hat{\beta}_{k}=0~\textrm{and}~\beta_{k}^{*}\neq0)$.

In below, we will use the conclusion of Lemma {\ref{lemma 2}} (KKT condition). Recall that $E[h_{k}(Z_{1})]=\sum_{j\in I_{*}}\beta_{j}^{*}<h_{k},h_{j}>=\sum_{j=1}^{W}\beta_{j}^{*}<h_{k},h_{j}>$. Since we assume that the density of $Z_{1}$ is the mixture density $h_{\beta^{*}}=\sum_{j\in I_{*}}\beta_{j}^{*}h_{j}$. So for $k\in I_{*}$ we have
\begin{align}\nonumber
&~~~~P(\hat{\beta}_{k}=0~\textrm{and}~\beta_{k}^{*}\neq0)=P\left(\left|\frac{1}{n}\sum_{i=1}^{n}h_{k}(X_{i})-\sum_{j=1}^{W}\hat{\beta}_{j}<h_{j},h_{k}>\right|\leq2\sqrt{2}v(\delta/2W)L_{k};\beta_{k}^{*}\neq0\right)\\\nonumber
&=P\left(\left|\frac{1}{n}\sum_{i=1}^{n}h_{k}(X_{i})-E[h_{k}(Z_{1})]+E[h_{k}(Z_{1})]-\sum_{j=1}^{W}\hat{\beta}_{j}<h_{j},h_{k}>\right|\leq2\sqrt{2}v(\delta/2W)L_{k};\beta_{k}^{*}\neq0\right)\\\nonumber
&=P\left(\left|\frac{1}{n}\sum_{i=1}^{n}h_{k}(X_{i})-E[h_{k}(Z_{1})]-\sum_{j=1}^{W}(\hat{\beta}_{j}-\beta_{j}^{*})<h_{j},h_{k}>\right|\leq2\sqrt{2}v(\delta/2W)L_{k};\beta_{k}^{*}\neq0\right)\\\nonumber
&= P\left(\left|\frac{1}{n}\sum_{i=1}^{n}h_{k}(X_{i})-E[h_{k}(Z_{1})]-\sum_{j\neq k}^{W}(\hat{\beta}_{j}
-\beta_{j}^{*})<h_{j},h_{k}>+\beta_{k}^{*}\|h_{k}\|^{2}\right|\leq2\sqrt{2}v(\delta/2W)L_{k}\right)\\\nonumber
&\leq P\left(|\beta_{k}^{*}\|h_{k}\|^{2} -2\sqrt{2}v(\delta/2W)L_{k}\leq\left|\frac{1}{n}\sum_{i=1}^{n}h_{k}(X_{i})-E[h_{k}(Z_{1})]\right|+\left|\sum_{j\neq k}^{W}(\hat{\beta}_{j}
-\beta_{j}^{*})<h_{j},h_{k}>\right|\right)\\\label{eq:11}
&\leq P\left(\left|\frac{1}{n}\sum_{i=1}^{n}h_{k}(X_{i})-E[h_{k}(Z_{1})]\right|\geq\frac{|\beta_{k}^{*}|\|h_{k}\|^{2}}{2}-\sqrt{2}v(\delta/2W)L_{k}\right)\\\label{eq:12}
&~~+P\left(\left|\sum_{j\neq k}^{W}(\hat{\beta}_{j}-\beta_{j}^{*})<h_{j},h_{k}>\right|\geq\frac{|\beta_{k}^{*}|\|h_{k}\|^{2}}{2}-\sqrt{2}v(\delta/2W)L_{k}\right).
\end{align}
Similar to Lemma {\ref{lemma 2}}, for (\ref{eq:11}), we use Hoeffding's inequality. Since $\|h_{k}\|=1$ for all $k$. Consider condition (B), $\min_{k\in I_{*}}|\beta^{*}_{k}|\geq4\sqrt{2}v(\delta/2W)L$ and $L\geq\max_{1\leq k\leq W}L_{k}$, then we have
\begin{align}\label{eq:13}
P&\left(\left|\frac{1}{n}\sum_{i=1}^{n}h_{k}(X_{i})-E[h_{k}(Z_{1})]\right|\geq\frac{|\beta_{k}^{*}|\|h_{k}\|^{2}}{2}-\sqrt{2}v(\delta/2W)L_{k}\right)\nonumber\\
&=P\left(\left|\frac{1}{n}\sum_{i=1}^{n}h_{k}(X_{i})-E[h_{k}(X_{1})]+E[h_{k}(X_{1})]-E[h_{k}(Z_{1})]\right|\geq\frac{|\beta_{k}^{*}|\|h_{k}\|^{2}}{2}-\sqrt{2}v(\delta/2W)L_{k}\right)\nonumber\\
&\leq P\left(\left|\frac{1}{n}\sum_{i=1}^{n}h_{k}(X_{i})-E[h_{k}(X_{1})]\right|\geq\frac{|\beta_{k}^{*}|}{2}-\sqrt{2}v(\delta/2W)L-\epsilon_{k}\right)~(\textrm{put}~\epsilon_{k}:=|E[h_{k}(X_{1})]-E[h_{k}(Z_{1})]|)\nonumber\\
&\leq P\left(\left|\frac{1}{n}\sum_{i=1}^{n}h_{k}(X_{i})-E[h_{k}(X_{1})]\right|\geq2\sqrt{2}v(\delta/2W)L-\sqrt{2}v(\delta/2W)L-\epsilon_{k}\right)\nonumber\\
&=P\left(\left|\frac{1}{n}\sum_{i=1}^{n}h_{k}(X_{i})-E[h_{k}(X_{1})]\right|\geq\sqrt{2}v(\delta/2W)L-\epsilon_{k}\right)\nonumber\\
&=P\left(\left|\frac{1}{n}\sum_{i=1}^{n}h_{k}(X_{i})-E[h_{k}(X_{1})]\right|\geq\sqrt{2}v(\delta/2W)L(1-\epsilon_{k}^{*})\right)~~(\textrm{let}~\epsilon_{k}^{*}=\epsilon_{k}/\sqrt{2}v(\delta/2W)L)\nonumber\\
&\leq2\exp\left\{-\frac{4n^{2}v^{2}(\delta/2W)L^{2}(1-\epsilon_{k}^{*})^2}{4nL^{2}}\right\}\nonumber\\
&=2\exp\left\{-n(1-\epsilon_{k}^{*})^2\frac{\log(2W^{2}/\delta)}{n}\right\}=2\left(\frac{\delta}{2W^{2}}\right)^{(1-\epsilon_{k}^{*})^2}.
\end{align}
For the upper bound of (\ref{eq:12}), using condition (A) and condition (B), by the definitions of $\rho_{*}(\beta^{*})$ and $W(\beta^{*})$, we have
\begin{align*}
&~~~~P\left(\left|\sum_{j\neq k}^{W}(\hat{\beta}_{j}-\beta_{j}^{*})<h_{j},h_{k}>\right|\geq\frac{|\beta_{k}^{*}|\|h_{k}\|^{2}}{2}-\sqrt{2}v(\delta/2W)L_{k}\right)\\
&=P\left(\left|\sum_{j\neq k}^{W}(\hat{\beta}_{j}-\beta_{j}^{*})<h_{j},h_{k}>\right|\geq\frac{|\beta_{k}^{*}|}{2}-\sqrt{2}v(\delta/2W)L_{k}\right)\\
&\leq P\left(\left|\sum_{j\neq k}^{W}(\hat{\beta}_{j}-\beta_{j}^{*})\frac{<h_{j},h_{k}>}{\|h_{j}\|\|h_{k}\|}\cdot\|h_{j}\|\|h_{k}\|\right|\geq2\sqrt{2}v(\delta/2W)L-\sqrt{2}v(\delta/2W)L\right)\\
&\leq P\left( {{\rho _*}({\beta ^*})\sum\limits_{j \ne k}^W {\left| {{{\hat \beta }_j} - \beta _j^*} \right|}  \ge \sqrt 2 v(\delta /2W)L} \right)\leq P\left( {\sum\limits_{j{\rm{ = 1}}}^W {\left| {{{\hat \beta }_j} - \beta _j^*} \right|}  \ge \frac{{\sqrt 2 v(\delta /2W)L}}{{{\rho _*}({\beta ^*})}}} \right)\\
&\leq P\left( {\sum\limits_{j{\rm{ = 1}}}^W {\left| {{{\hat \beta }_j} - \beta _j^*} \right|}  \ge \frac{{{\rm{288}}\sqrt 2 {G^*}v(\delta /2W)}}{{{L_{\min }}{\lambda _W}}}} \right)\leq \frac{\delta}{W}.
\end{align*}
Where the second last inequality is by condition (A) and the last inequality above is by using the $\ell_1$-estimation oracle inequality in Corollary {\ref{corollary 2}}.

Therefore, by the definition of $W(\beta^{*})$, $W(\beta^{*})=|I_{*}|\leq W$, we have
\begin{align*}
P(I_{*}\nsubseteqq\hat{I})\leq W(\beta^{*})\max_{k\in I_{*}}P(\hat{\beta}_{k}=0)&\leq W(\beta^{*})2\left(\frac{\delta}{2W^{2}}\right)^{(1-\epsilon_{k}^{*})^2}+W(\beta^{*})\frac{\delta}{W}\\
\leq & 2W\left(\frac{\delta}{2W^{2}}\right)^{(1-\epsilon_{k}^{*})^2}+W\frac{\delta}{W}=2W\left(\frac{\delta}{2W^{2}}\right)^{(1-\epsilon_{k}^{*})^2}+\delta.
\end{align*}

For the control of $P(\hat{I}\nsubseteqq I_{*})$, let
\begin{align}\label{eq:15}
\tilde{\eta}=\underset{\eta\in \mathbb{R}^{W(\beta^{*})}}{\arg\min}z(\eta),
\end{align}
where
\begin{align*}
z(\eta)=-\frac{2}{n}\sum_{i=1}^{n}\sum_{j\in I_{*}}\eta_{j}h_{j}(X_{i})+\|\sum_{j\in I_{*}}\eta_{j}h_{j}\|^{2}+\sum_{j\in I_{*}}(4\sqrt{2}v(\delta /2)L_{j}+2cB)|\eta_{j}|+c\sum_{j\in I_{*}}\eta_{j}^{2}.
\end{align*}

Consider the following random event,
\begin{align}\label{eq:16}\nonumber
\bigcap_{k\notin I_{*}}&\left\{\left|-\frac{1}{n}\sum_{i=1}^{n}h_{k}(X_{i})+\sum_{j\in I_{*}}\tilde{\eta}_{j}<h_{j},h_{k}>\right|\leq2\sqrt{2}v(\delta /2)L_{k}\right\}\\
&\subseteq\bigcap_{k\notin I_{*}}\left\{\left|-\frac{1}{n}\sum_{i=1}^{n}h_{k}(X_{i})+\sum_{j\in I_{*}}\tilde{\eta}_{j}<h_{j},h_{k}>\right|\leq2\sqrt{2}v(\delta /2W)L\right\}:=\Psi.
\end{align}
Let $\bar{\eta}\in \mathbb{R}^{W}$ be a vector corresponding to the component of the index set $I_{*}$ having $\tilde{\eta}$ given by equation (\ref{eq:15}), and the component at other corresponding positions is 0. By Lemma {\ref{lemma 1}}, we know that $\bar{\eta}\in \mathbb{R}^{W}$ is a solution of (\ref{eq:9}) on the event $\Psi$. It is recalled that $\hat{\beta}\in \mathbb{R}^{W}$ is also a solution of (\ref{eq:9}). Through the definition of the indicator set $\hat{I}$, we have $\hat{\beta}_{k}\neq0$ for $k\in \hat{I}$. By construction, we obtain $\tilde{\eta}_{k}\neq0$ for some subset $T\subseteqq I_{*}$. The KKT conditions indicate that any two solutions have non-zero components at the same positions. Therefore, $\hat{I}=T\subseteqq I_{*}$ on the event $\Psi$. Further, we have
\begin{align}\nonumber
&~~~~P(\hat{I}\nsubseteq I_{*})\leq P(\Psi^{c})
=P\left(\bigcup_{k\notin I_{*}}\left\{\left|-\frac{1}{n}\sum_{i=1}^{n}h_{k}(X_{i})+\sum_{j\in I_{*}}\tilde{\eta}_{j}<h_{j},h_{k}>\right|\geq2\sqrt{2}v(\delta /2W)L\right\}\right)\\\nonumber
&\leq\sum_{k\notin I_{*}}P\left\{\left|-\frac{1}{n}\sum_{i=1}^{n}h_{k}(X_{i})+\sum_{j\in I_{*}}\tilde{\eta}_{j}<h_{j},h_{k}>\right|\geq2\sqrt{2}v(\delta /2W)L\right\}\\\nonumber
&=\sum_{k\notin I_{*}}P\left\{\left|-\frac{1}{n}\sum_{i=1}^{n}h_{k}(X_{i})+E[h_{k}(Z_{1})]-E[h_{k}(Z_{1})]+\sum_{j\in I_{*}}\tilde{\eta}_{j}<h_{j},h_{k}>\right|\geq2\sqrt{2}v(\delta /2W)L\right\}\\\nonumber
&=\sum_{k\notin I_{*}}P\left\{\left|\frac{1}{n}\sum_{i=1}^{n}h_{k}(X_{i})-E[h_{k}(Z_{1})]-\sum_{j\in I_{*}}(\tilde{\eta}_{j}-\beta_{j}^{*})<h_{j},h_{k}>\right|\geq2\sqrt{2}v(\delta /2W)L\right\}\\\label{eq:17}
&\leq\sum_{k\notin I_{*}}P\left\{\left|\frac{1}{n}\sum_{i=1}^{n}h_{k}(X_{i})-E[h_{k}(Z_{1})]\right|\geq\sqrt{2}v(\delta /2W)L\right\}\\\label{eq:18}
&~~+\sum_{k\notin I_{*}}P\left\{\sum_{j\in I_{*}}|\tilde{\eta}_{j}-\beta_{j}^{*}|\left|<h_{j},h_{k}>\right|\geq\sqrt{2}v(\delta /2W)L\right\}.
\end{align}
According to the previously proved (\ref{eq:13}) formula, we have
\begin{align*}
&~~~~\sum_{k\notin I_{*}}P\left\{\left|\frac{1}{n}\sum_{i=1}^{n}h_{k}(Z_{i})-Eh_{k}(Z_{1})\right|\geq\sqrt{2}v(\delta /2W)L\right\}\\
&\leq\sum_{k\notin I_{*}}P\left\{\left|\frac{1}{n}\sum_{i=1}^{n}h_{k}(X_{i})-Eh_{k}(X_{1})\right|\geq\sqrt{2}v(\delta /2W)L-\epsilon_{k}\right\}\\
&=\sum_{k=1}^{W}P\left\{\left|\frac{1}{n}\sum_{i=1}^{n}h_{k}(X_{i})-Eh_{k}(X_{1})\right|\geq\sqrt{2}v(\delta /2W)L(1-\epsilon_{k}^{*})\right\}\leq 2W(\frac{\delta}{2W^{2}})^{(1-\epsilon_{k}^{*})^{2}}.
\end{align*}
For the upper bound of (\ref{eq:18}), observe Theorem {\ref{theorem 2}}, we can use a larger $v(\delta/2W)$ instead of $v(\delta/2)$. Considering the construction of $\tilde{\eta}$ in (\ref{eq:15}), we find
\begin{align*}
P\left(\sum_{j\in I_{*}}|\tilde{\eta}_{j}-\beta_{j}^{*}|\geq\frac{{{\rm{288}}\sqrt 2 {G^*}v(\delta /2W)}}{{{L_{\min }}{\lambda _W}}}\right)\leq \frac{\delta}{W}.
\end{align*}
Similarly, we have
\begin{align*}
&~~~~\sum_{k\notin I_{*}}P\left\{\sum_{j\in I_{*}}|\tilde{\eta}_{j}-\beta_{j}^{*}|\left|<h_{j},h_{k}>\right|\geq\sqrt{2}v(\delta /2W)L\right\}\\
&\leq\sum_{k=1}^{W}P\left\{\sum_{j\in I_{*}}|\tilde{\eta}_{j}-\beta_{j}^{*}|\left|\frac{<h_{j},h_{k}>}{\|h_{j}\|\|h_{k}\|}\|h_{j}\|\|h_{k}\|\right|\geq\sqrt{2}v(\delta /2W)L\right\}\\
&\leq\sum_{k=1}^{W}P\left\{\sum_{j\in I_{*}}|\tilde{\eta}_{j}-\beta_{j}^{*}|{{{\rho _*}({\beta ^*})}}\geq\sqrt{2}v(\delta /2W)L\right\}\\
&=\sum_{k=1}^{W}P\left\{\sum_{j\in I_{*}}|\tilde{\eta}_{j}-\beta_{j}^{*}|\geq \frac{{\sqrt 2 v(\delta /2W)L}}{{{\rho _*}({\beta ^*})}}\right\}\\
(\text{using condition (A)})~&\leq\sum_{k=1}^{W}P\left\{\sum_{j\in I_{*}}|\tilde{\eta}_{j}-\beta_{j}^{*}|\geq\frac{{{\rm{288}}\sqrt 2 {G^*}v(\delta /2W)}}{{{L_{\min }}{\lambda _W}}}\right\}\leq \sum_{k=1}^{W}\frac{\delta}{W}=\delta.
\end{align*}
In summary, we can get
\begin{align*}
P(\hat{I}\neq I_{*})\leq P(I_{*}\nsubseteqq\hat{I})+P(\hat{I}\nsubseteqq I_{*})&\leq 2W\left(\frac{\delta}{2W^{2}}\right)^{(1-\epsilon_{k}^{*})^2}+\delta+2W\left(\frac{\delta}{2W^{2}}\right)^{(1-\epsilon_{k}^{*})^2}+\delta\\
&=4W\left(\frac{\delta}{2W^{2}}\right)^{(1-\epsilon_{k}^{*})^2}+2\delta.
\end{align*}
This completes the proof of Theorem {\ref{theorem 3}}.

\section*{Acknowledgements}
The authors would like to thank Song Xi Chen's Group (\url{http://songxichen.gsm.pku.edu.cn/}) for sharing the meteorological data sets.


\begin{thebibliography}{999}
\providecommand{\natexlab}[1]{#1} \providecommand{\url}[1]{\texttt{#1}} \providecommand{\urlprefix}{URL }
\expandafter\ifx\csname urlstyle\endcsname\relax
  \providecommand{\doi}[1]{doi:\discretionary{}{}{}#1}\else
  \providecommand{\doi}{doi:\discretionary{}{}{}\begingroup
  \urlstyle{rm}\Url}\fi
\providecommand{\eprint}[2][]{\url{#2}}
\small{
\bibitem[{Balakrishnan et al.(2017)}]{Balakrishnanl17}
Balakrishnan, S., Wainwright, M. J., \& Yu, B. (2017). Statistical guarantees for the EM algorithm: From population to sample-based analysis. The Annals of Statistics, 45(1), 77-120.

\bibitem[{Belloni et al.(2017)}]{Belloni17}
Belloni, A., Rosenbaum, M., \& Tsybakov, A. B. (2017). Linear and conic programming estimators in high dimensional errors-in-variables models. Journal of the Royal Statistical Society: Series B (Statistical Methodology), 79(3), 939-956.

\bibitem[{Bertin et al.(2011)}]{Bertin11}
Bertin, K., Le Pennec, E., \& Rivoirard, V. (2011). Adaptive Dantzig density estimation. Annales de l'IHP Probabilit{\'e}s et statistiques (Vol. 47, No. 1, pp. 43-74).

\bibitem[Biau and Devroye(2005)]{Biau2005}
Biau, G., \& Devroye, L. (2005). Density estimation by the penalized combinatorial method. Journal of Multivariate Analysis, 94(1), 196-208.

\bibitem[{Bickel et al.(2009)}]{Bickel09}
Bickel, P. J., Ritov, Y. A., \& Tsybakov, A. B. (2009). Simultaneous analysis of Lasso and Dantzig selector. The Annals of Statistics, 1705-1732.

\bibitem[Buhlmann and van de Geer(2011)]{Buhlmann2011}
Buhlmann, P., \& van de Geer, S. (2011). Statistics for high-dimensional data: methods, theory and applications. Springer.

\bibitem[{Boucheron et al.(2013)}]{Boucheron13}
Boucheron, S., Lugosi, G., \& Massart, P. (2013). Concentration inequalities: A nonasymptotic theory of independence. Oxford university press.

\bibitem[{Bunea(2008)}]{Bunea08}
Bunea, F. (2008). Honest variable selection in linear and logistic regression models via $\ell_{1}$ and $\ell_{1}+\ell_{2}$ penalization. Electronic Journal of Statistics, 2, 1153-1194.

\bibitem[{Bunea et al.(2010)}]{Bunea10}
Bunea, F., Tsybakov, A. B., Wegkamp, M. H., \& Barbu, A. (2010). Spades and mixture models. The Annals of Statistics, 38(4), 2525-2558.

\bibitem[{Candes(2008)}]{Candes08}
Candes, E. J. (2006). Modern statistical estimation via oracle inequalities. Acta numerica, 15, 257-325.


\bibitem[{Cheng and van Ness(1999)}]{Cheng1999}
Cheng, C. L. \& van Ness, J. W. (1999) Statistical Regression with Measurement Error. New York: Wiley.

\bibitem[{Chen(1998)}]{Chen1998}
Chen, S. X. (1998). Measurement errors in line transect surveys. Biometrics, 899-908.

\bibitem[{Chen and Khalili(2008)}]{Chen2008}
Chen, J., \& Khalili, A. (2008). Order selection in finite mixture models with a nonsmooth penalty. Journal of the American Statistical Association, 103(484), 1674-1683.

\bibitem[{Chow and Teicher(2003)}]{chow1997probability}
Chow, Y. S., \& Teicher, H. (2003). Probability theory: independence, interchangeability, martingales, 3rd. Springer.


\bibitem[{DasGupta(2008)}]{Das08}
DasGupta, A. (2008). Asymptotic theory of statistics and probability. Springer.

\bibitem[{Devroye and Lugosi(2001)}]{Lugosi2001}
Devroye, L., \& Lugosi, G. (2001). Combinatorial methods in density estimation. Springer.

\bibitem[{Donoho and Johnstone(1994)}]{Donoho1994}
Donoho, D. L., \& Johnstone, J. M. (1994). Ideal spatial adaptation by wavelet shrinkage. biometrika, 81(3), 425-455.

\bibitem[{Hall et al.(2008)}]{Hall08}
Hall, P., \& Lahiri, S. N. (2008). Estimation of distributions, moments and quantiles in deconvolution problems. The Annals of Statistics, 36(5), 2110-2134.

\bibitem[{Hersbach et al. (2018)}]{Hersbach2018}
Hersbach, H. et al. (2018). Operational global reanalysis: progress, future directions and synergies with NWP[M]. European Centre for Medium Range Weather Forecasts.

\bibitem[{Lemler(2016)}]{Lemler2016}
Lemler, S. (2016). Oracle inequalities for the Lasso in the high-dimensional Aalen multiplicative intensity model. Annales De L Institut Henri Poincare-probabilites Et Statistiques, 52(2), 981-1008.

\bibitem[Meister(2009)]{Martin09}
Martin, R. (2009). Fast Nonparametric Estimation of a Mixing Distribution with Application
to High Dimensional Inference, PhD thesis, Purdue University.

\bibitem[{McLachlan et al.(2019)}]{McLachlan19}
McLachlan, G. J., Lee, S. X., \& Rathnayake, S. I. (2019). Finite mixture models. Annual review of statistics and its application, 6, 355-378.

\bibitem[Meister(2006)]{Meister2006}
Meister, A. (2006). Density estimation with normal measurement error with unknown variance. Statistica Sinica, 195-211.

\bibitem[Nakamura(1990)]{Nakamura90}
Nakamura, T. (1990). Corrected score function for errors-in-variables models: Methodology and application to generalized linear models. Biometrika, 77(1), 127-137.

\bibitem[{Rosenbaum and Tsybakov(2010)}]{Rosenbaum2010}
Rosenbaum, M., \& Tsybakov, A. B. (2010). Sparse recovery under matrix uncertainty. The Annals of Statistics, 2620-2651.

\bibitem[{Schennach and Bonhomme(2013)}]{Schennach2013}
Schennach, S., \& Bonhomme, S. (2013). Penalized Least Squares Methods for Latent Variables Models. In Advances in Economics and Econometrics: Volume 3, Econometrics: Tenth World Congress (Vol. 51, p. 338). Cambridge University Press.

\bibitem[{Tokdar et al.(2009)}]{Tokdar09}
Tokdar, S. T., Martin, R., \& Ghosh, J. K. (2009). Consistency of a recursive estimate of mixing distributions. The Annals of Statistics, 37(5A), 2502-2522.

\bibitem[Tibshirani(1996)]{Tibshirani1996}
Tibshirani, R. (1996). Regression shrinkage and selection via the lasso. Journal of the Royal Statistical Society. Series B (Methodological), 267-288.

\bibitem[{Zhang and Jia(2020)}]{Zhang20}
Zhang H., \& Jia, J. (2020). Elastic-net Regularized High-dimensional Negative Binomial Regression: Consistency and Weak Signals Detection. Statistica Sinica, \url{https://doi.org/10.5705/ss.202019.0315}


\bibitem[{Zou and Zhang(2005)}]{Zou2005}
Zou, H., \& Hastie, T. (2005). Regularization and variable selection via the elastic net. Journal of the Royal Statistical Society: Series B (Statistical Methodology), 67(2), 301-320.

\bibitem[{Zou and Zhang(2009)}]{Zou2009}
Zou, H., \& Zhang, H. (2009). On the adaptive elastic-net with a diverging number of parameters. Annals of statistics, 37(4), 1733.


}
\end{thebibliography}
\end{document}